**RESEARCH ARTICLE**

# Motohashi's formula for the fourth moment of individual Dirichlet *L*-functions and applications

Ikuya Kaneko

The Division of Physics, Mathematics and Astronomy, California Institute of Technology, 1200 E. California Blvd., Pasadena, CA, 91125, USA; E-mail: ikuyak@icloud.com.

**Abstract**
A new reciprocity formula for Dirichlet *L*-functions associated to an arbitrary primitive Dirichlet character of prime modulus $q$ is established. We find an identity relating the fourth moment of individual Dirichlet *L*-functions in the *t*-aspect to the cubic moment of central *L*-values of Hecke–Maaß newforms of level at most $q^2$ and primitive central character $\psi^2$ averaged over all primitive nonquadratic characters $\psi$ modulo $q$. Our formula can be thought of as a reverse version of recent work of Petrow–Young. Direct corollaries involve a variant of Iwaniec's short interval fourth moment bound and the twelfth moment bound for Dirichlet *L*-functions, which generalise work of Jutila and Heath-Brown, respectively. This work traverses an intersection of classical analytic number theory and automorphic forms.

## Contents

## 1. Introduction

Estimating moments of families of *L*-functions is a central problem in analytic number theory not only due to their substantial applications, but also since they give an insight into the behaviour of *L*-functions in the critical strip. Our interests lie in the $2k$th moments of the Riemann zeta and Dirichlet *L*-functions:

$$\mathcal{Z}_k(g) = \int_{-\infty}^{\infty} \left| \zeta\left(\frac{1}{2} + it\right) \right|^{2k} g(t)\,dt, \qquad \mathcal{Z}_k(g;\chi) = \int_{-\infty}^{\infty} \left| L\left(\frac{1}{2} + it, \chi\right) \right|^{2k} g(t)\,dt,$$

where $k \geqslant 1$ is an integer and the test function *g* is of rapid decay. The initial cases $k = 1, 2$ have been successfully studied for $\mathcal{Z}_k(g)$, and the other cases have remained untouched thus far. In this article, we establish Motohashi's formula (also known as a spectral reciprocity formula) for the fourth moment of Dirichlet *L*-functions, which was unsolved since being posed as a question of study by Motohashi in 1992.

### 1.1. Overview and motivation

In the 1990s, Motohashi [55, Theorem 4.2] pondered a mysterious identity relating the smoothed fourth moment of the Riemann zeta function $\zeta(1/2+it)$ to a spectral cubic moment of automorphic *L*-functions associated to the group $\mathrm{SL}_2(\mathbb{Z})$. We assume a fixed test function *g* to be of Schwartz class, and we denote by $\mathcal{B}(q, \chi)$ the set of Hecke–Maaß forms of level *q* and central character $\chi$; we write $\mathcal{B}(\Gamma_0(q))$ as usual when $\chi$ is principal. His formula then asserts that the following spectral decomposition holds up to an explicit description of holomorphic, Eisenstein and residual contributions that we shall elide here:

$$\int_{-\infty}^{\infty} \left| \zeta\left(\frac{1}{2} + it\right) \right|^4 g(t)\,dt \leftrightsquigarrow \sum_{f \in \mathcal{B}(\Gamma_0(1))} L\left(\frac{1}{2}, f\right)^3 \check{g}(t_f), \tag{1.1}$$

where $\check{g}$ is an elaborate integral transform of *g* involving the Gauß hypergeometric function. The right-hand side of equation (1.1) must be understood as a complete integral over the full spectrum of level 1 automorphic forms, including holomorphic, discrete and continuous spectra. Motohashi has given several approaches in the spirit of analytic number theory and representation theory. Note that all of his methods are in the framework of relative trace formulæ. On the other hand, Michel–Venkatesh ([49, §4.3.3], [50, §4.5.3]) suggested an elegant geometric and spectral stratagem to substantiate equation (1.1). Following their perspective, Nelson [58] studied the cubic moment of automorphic *L*-functions on $\mathrm{PGL}_2$ via the use of the regularised diagonal periods of products of Eisenstein series. We also refer the reader to work of Wu et al. [4, 68].

In terms of automorphic representation theoretic language, the spectral reciprocity formula in equation (1.1) is a connection between the fourth moment of $\mathrm{GL}_1$ $L$-functions (geometric side) and the cubic moment of $\mathrm{GL}_2$ $L$-functions (spectral side). These sides are derived from an application of the relative trace formula. The advancement toward Motohashi's formula is summarised in the following table.[1]

| $L$-functions | Individual | Character average |
| --- | --- | --- |
| Riemann zeta function | Bruggeman–Motohashi [19, 53, 55] | n/a |
| Dirichlet $L$-functions | Theorem 1.1 | Blomer et al. [10] |
| Dedekind zeta functions | Bruggeman–Motohashi [17, 18, 54] | unknown |

These works share similar structure to Motohashi's formula in equation (1.1) (compare [56]). In particular, the heuristics of Blomer et al. [10] establishes that the character average of the smoothed fourth moment of Dirichlet $L$-functions weighted by $\chi(a)\overline{\chi}(b)$ is expressed by means of a cubic moment of $L$-functions associated to Hecke–Maaß newforms of level $ab$. Their method relies on brute force calculations and is cleverer than in [55] in the sense that their reasoning is rather symmetric and utilises additive reciprocity. There exist several articles in the antecedent literature on such a fourth moment problem; see [30, 41, 69]. Although these achievements address the problem of obtaining an asymptotic, one must go through a harder route to prove Motohashi's formula.

There are two versions of spectral reciprocity formulæ: $\mathrm{GL}_4 \times \mathrm{GL}_2 \leftrightsquigarrow \mathrm{GL}_4 \times \mathrm{GL}_2$ and $\mathrm{GL}_2 \times \mathrm{GL}_2 \leftrightsquigarrow \mathrm{GL}_3 \times \mathrm{GL}_2$. The former instance involves work of Andersen–Kıral [1], Blomer–Khan [11, 12], Blomer–Li–Miller [13], Humphries–Khan [33], Jana–Nunes [40], Kuznetsov [47], Nunes [59] and Zacharias [71]. The latter involves work of Blomer et al. [10], Nelson [58], Petrow [62], Petrow–Young [63, 64], Wu [68] and Young [69]. It behooves one to touch on the article in preparation due to Humphries–Khan to deduce $\mathrm{GL}_3 \times \mathrm{GL}_2 \leftrightsquigarrow \mathrm{GL}_4 \times \mathrm{GL}_1$ reciprocity. This renders an extension of Motohashi's formula in equation (1.1).

## 1.2. *Statement of main result*

Let $\chi$ be an arbitrary primitive Dirichlet character modulo an integer $q$, and let $\mathcal{R}_4^+$ be a subdomain of $\mathbb{C}^4$, where all four parameters have real parts greater than one. In this article, we update the progress toward spectral reciprocity formulæ and extend equation (1.1) to the fourth moment of *individual* Dirichlet $L$-functions associated to $\chi$ in the $t$-aspect. Thus the setup that we build upon is as follows. For $q \in \mathbb{N}$, we consider

$$\mathcal{Z}_2(g; \chi) = \int_{-\infty}^{\infty} \left| L\left(\frac{1}{2} + it, \chi\right) \right|^4 g(t) dt. \tag{1.2}$$

This is seen as a character analogue of equation (1.1). The twist by $\chi$ substantially complicates our analysis, and we will encounter various intricate character sums. If $g$ is a sufficiently nice test function, then we define

$$\mathcal{Z}_2(s_1, s_2, s_3, s_4; g; \chi) = \int_{-\infty}^{\infty} L(s_1 + it, \chi) L(s_2 + it, \chi) L(s_3 - it, \overline{\chi}) L(s_4 - it, \overline{\chi}) g(t) dt, \tag{1.3}$$

where the parameters or shifts satisfy $(s_1, s_2, s_3, s_4) \in \mathcal{R}_4^+$. In order to establish Motohashi's formula for the fourth moment in equation (1.2), we initially work with equation (1.3) by exploiting Dirichlet series expansions for the integrand in the region of absolute convergence, and we take the limit $(s_1, s_2, s_3, s_4) \to (1/2, 1/2, 1/2, 1/2)$ after the meromorphic continuation as in the original work of Motohashi [55].

---

[1] The character average for the fourth moment of Dedekind zeta functions should be meant only for quadratic number fields. In this scenario, the problem could be solved with classical technology, although this has not been worked out anywhere.

We proceed to the rigorous statement of our reciprocity formula, which requires a bit of notation. For simplicity, we assume that $q$ is prime. The classification of automorphic forms is convenient in the sequel:

○ Cuspidal holomorphic newforms $f$ of weight $k \equiv \kappa(\psi) = \kappa \pmod 2$, level $q$, central character $\psi$ and Hecke eigenvalues $\lambda_f(n) \in \mathbb{C}$; we denote the set of such forms by $\mathcal{B}_k^*(q, \psi)$;

○ Cuspidal Maaß newforms $f$ of spectral parameter $t_f \in \mathbb{R} \cup [-i\vartheta, i\vartheta]$, weight $\kappa = (1 - \psi(-1))/2$, level $q$, central character $\psi$ and Hecke eigenvalues $\lambda_f(n) \in \mathbb{C}$, where at the current state of knowledge, $\vartheta = 7/64$ can be taken (see [45]; although $\vartheta = 0$ is expected); we denote the set of such forms by $\mathcal{B}_\kappa^*(q, \psi)$;

○ Unitary Eisenstein series $E(z, s, f)$, where $s = 1/2 + it$ with $t \in \mathbb{R} \setminus \{0\}$ and $\mathcal{B}(\psi_1, \psi_2) \ni f$ with $\psi = \psi_1 \psi_2$ is a certain finite set depending upon $\psi_1, \psi_2$ corresponding to an orthonormal basis in the space of the induced representation constructed out of $(\psi_1, \psi_2)$. Their $n$th Hecke eigenvalue is written as $\lambda_f(n, t) = \sum_{ab=n} \psi_1(a) a^{it} \psi_2(b) b^{-it}$ for $(n, q) = 1$.

Let $\tau(\chi)$ and $J(\chi, \psi)$ be the Gauß sum and the Jacobi sum, respectively, so that we define $\mathcal{H}_\pm(\chi, \psi) = \psi(\mp 1) \tau(\psi) \tau(\overline{\chi \psi}) J(\chi, \psi)$. For a test function as in Convention 1.10, we define the spectral mean values[2]

$$\mathcal{J}_\pm^{\text{Maaß}} := \frac{q^{s_2 - s_1 - 2}}{\tau(\overline{\chi})} \sum_{d \mid q} \frac{\mu(d)}{d} \sum_{\psi \pmod q}^{\#} \mathcal{H}_\pm(\chi, \psi) \sum_{f \in \mathcal{B}_\kappa^*(dq, \psi^2)} \epsilon_f^{(1 \mp 1)/2}$$

$$\times \frac{L\left(\frac{1 - s_1 + s_2 + s_3 - s_4}{2}, f \otimes \overline{\psi}\right) L\left(\frac{1 - s_1 + s_2 - s_3 - s_4}{2}, f \otimes \overline{\psi}\right) L\left(\frac{s_1 + s_2 + s_3 + s_4 - 1}{2}, f \otimes \overline{\psi}\right)}{L(1, \text{Ad}^2 f)} \Phi_{\mathsf{s}}^\pm(it_f),$$

$$\mathcal{J}_\pm^{\text{Eis}} := \frac{2q^{s_2 - s_1 - 1}}{\tau(\overline{\chi}) \varphi(q)} \sum_{d \mid q} \frac{\mu(d)}{d} \sum_{\psi \pmod q}^{\#} \mathcal{H}_\pm(\chi, \psi)$$

$$\times \sum_{\psi_1 \psi_2 = \psi^2} \sum_{f \in \mathcal{B}(\psi_1, \psi_2)} \int_{-\infty}^\infty \frac{\mathcal{S}_f(t; s_1, s_2, s_3, s_4)}{|L(1 + 2it, \overline{\psi_1} \psi_2)|^2} \Phi_{\mathsf{s}}^\pm(it) \frac{dt}{2\pi},$$

$$\mathcal{J}_+^{\text{hol}} := \frac{q^{-2}}{\tau(\overline{\chi})} \left(\frac{2\pi}{q}\right)^{s_1 - s_2} \cos\left(\frac{\pi(s_3 - s_4)}{2}\right) \sum_{d \mid q} \frac{\mu(d)}{d} \sum_{\psi \pmod q}^{\#} \mathcal{H}_+(\chi, \psi) \sum_{\substack{k > \kappa \\ k \equiv \kappa \pmod 2}} \sum_{f \in \mathcal{B}_k^*(dq, \psi^2)} i^k$$

$$\times \frac{L\left(\frac{1 - s_1 + s_2 + s_3 - s_4}{2}, f \otimes \overline{\psi}\right) L\left(\frac{1 - s_1 + s_2 - s_3 + s_4}{2}, f \otimes \overline{\psi}\right) L\left(\frac{s_1 + s_2 + s_3 + s_4 - 1}{2}, f \otimes \overline{\psi}\right)}{L(1, \text{Ad}^2 f)} \Xi_{\mathsf{s}}\left(\frac{k-1}{2}\right),$$

where $\mathcal{J}_\pm^{\text{Maaß}}$ and $\mathcal{J}_\pm^{\text{hol}}$ involve the cubic moment of twisted automorphic $L$-functions and the continuous term $\mathcal{J}_\pm^{\text{Eis}}$ should be regarded as the sixth moment of Dirichlet $L$-functions. Here, $\#$ on the sum signifies that the sum runs over all primitive nonquadratic characters modulo $q$, and we define

$$\mathcal{S}_f(t; s_1, s_2, s_3, s_4) = L\left(\frac{1 - s_1 + s_2 + s_3 - s_4}{2} + it, \overline{\psi} \psi_2\right) L\left(\frac{1 - s_1 + s_2 + s_3 - s_4}{2} - it, \overline{\psi} \psi_1\right)$$

$$\times L\left(\frac{1 - s_1 + s_2 - s_3 + s_4}{2} + it, \overline{\psi} \psi_2\right) L\left(\frac{1 - s_1 + s_2 - s_3 + s_4}{2} - it, \overline{\psi} \psi_1\right) \quad (1.4)$$

$$\times L\left(\frac{s_1 + s_2 + s_3 + s_4 - 1}{2} + it, \psi \overline{\psi_2}\right) L\left(\frac{s_1 + s_2 + s_3 + s_4 - 1}{2} - it, \psi \overline{\psi_1}\right).$$

---

[2]These are the main spectral contributions to the cubic moment side. One should take care of degenerate terms stemming from the principal and quadratic Dirichlet characters in the proof of Theorem 1.1, but they are in principle the same size as the main spectral contributions in terms of estimations. Hence one can disregard them when using Motohashi's formula to deduce certain moment bounds and asymptotic estimations.

The factor $\epsilon_f$ in the summand denotes the parity of a Maaß form $f$. In addition, $\Phi_{\mathsf{s}}^{\pm}$ and $\Xi_{\mathsf{s}}$ are introduced in equations (3.26), (3.27) and (3.28), respectively. The Dirichlet series expansions of the adjoint square $L$-function $L(s, \mathrm{Ad}^2 f)$, and the twisted automorphic $L$-function $L(s, f \otimes \psi)$ in $\Re(s) > 1$ (that serve as a definition of these functions) are given in equations (2.24) and (2.25), respectively. We are now ready to describe the reciprocity formula to which we have alluded previously.

**Theorem 1.1.** *Let $\chi$ be a primitive Dirichlet character modulo a prime $q$, and let $\mathsf{s} = (s_1, s_2, s_3, s_4)$ with $\overline{\mathsf{s}} = (\overline{s_3}, \overline{s_4}, \overline{s_2}, \overline{s_1})$. Let $(s_1, s_2, s_3, s_4) \in \mathbb{C}^4$ be such that $|s_1|, |s_2|, |s_3|, |s_4| < B$ with $B$ sufficiently large. If a test function $g$ satisfies the basic assumption in Convention 1.10, then we have that*

$$\mathcal{Z}_2(\mathsf{s}; g; \chi) = \mathcal{N}(\mathsf{s}; g; \chi) + \sum_{\pm} \left( \mathcal{J}_{\pm}(\mathsf{s}; g; \chi) + \mathcal{E}_{\pm}(\mathsf{s}; g; \chi) + \overline{\mathcal{J}_{\pm}(\overline{\mathsf{s}}; g; \chi)} + \overline{\mathcal{E}_{\pm}(\overline{\mathsf{s}}; g; \chi)} \right), \qquad (1.5)$$

*where $\mathcal{N}$ is an explicitly computable main term, and we decompose*

$$\mathcal{J}_{\pm}(\mathsf{s}; g; \chi) = \{\mathcal{J}_{\pm}^{\mathrm{Maa}} + \mathcal{J}_{\pm}^{\mathrm{Eis}} + \delta_{\pm=+} \mathcal{J}_{+}^{\mathrm{hol}}\}(\mathsf{s}; g; \chi),$$
$$\mathcal{E}_{\pm}(\mathsf{s}; g; \chi) = \{\mathcal{E}_{\pm}^{\mathrm{Maa}} + \mathcal{E}_{\pm}^{\mathrm{Eis}} + \delta_{\pm=+} \mathcal{E}_{+}^{\mathrm{hol}}\}(\mathsf{s}; g; \chi).$$

*Here $\mathcal{E}_{\pm}(\mathsf{s}; g; \chi)$ is the degenerate term defined in Section 3.5.2, which has similar shape to $\mathcal{J}_{\pm}(\mathsf{s}; g; \chi)$.*

It should be feasible to generalise Theorem 1.1 to the case of any positive integer $q$. If $q$ is not prime, we need to utilise a general version of the transformation formula due to Blomer–Milićević [15, (2.2)]:

$$\sum_{(c,q)=1} \chi(c) S(m, n; c) h(c) = \frac{\chi_1(m)}{\tau(\chi_1)} \sum_{d|q} \mu(d) \sum_{dq_1|c} S_{\chi_1}(m, nq_1^2; c) h\left(\frac{c}{q_1}\right). \qquad (1.6)$$

This is valid for an arbitrary Dirichlet character modulo $q$ induced from a primitive character $\chi_1$ modulo $q_1 \mid q$ and for every $(m, q_1) = 1$. The formula in equation (1.6) translates Kloosterman sums associated to the $(\infty, 0)$ cusp-pair into twisted Kloosterman sums associated to the $(\infty, \infty)$ cusp-pair. The bottleneck is that the Kloosterman sum on the right-hand side of equation (1.6) involves $q_1^2$ instead of $q^2$.

**Remark 1.2.** In the author's recent work [44], the second moment of the product of the Riemann zeta and Dirichlet $L$-functions was contemplated:

$$\int_{-\infty}^{\infty} \left| \zeta\left(\frac{1}{2} + it\right) L\left(\frac{1}{2} + it, \chi\right) \right|^2 g(t) dt. \qquad (1.7)$$

If one replaces the Riemann zeta function with a Dirichlet $L$-function, this is in accordance with equation (1.2). Although the definition in equation (1.2) is similar to equation (1.7), the cubic moment side in Theorem 1.1 has quite different shape. This is due to the occurrence of additional character sums when we consider the fourth moment of Dirichlet $L$-functions. However, the resulting form of an integral transform of $g$ never becomes altered.

### 1.3. *Quantitative applications*

We are able to provide various quantitative applications of Theorem 1.1. A hybrid fourth moment bound of interval $H$ for individual Dirichlet $L$-functions is initially proven.

**Corollary 1.3.** *Let $T^{1/2} \leqslant H \leqslant T(\log T)^{-1}$. For any primitive Dirichlet character $\chi$ modulo a prime $q$, we have*

$$\int_T^{T+H} \left| L\left(\frac{1}{2} + it, \chi\right) \right|^4 dt \ll_{\epsilon} H^{1+\epsilon} q^{\epsilon} + \left(\frac{qT}{\sqrt{H}}\right)^{1+\epsilon}. \qquad (1.8)$$

The proof of Corollary 1.3 looks circular at first glance: we use work of Petrow–Young [63] on cubic moments; they show that bounds for these cubic moments follow from bounds for the fourth moment of Dirichlet $L$-functions. Nonetheless, this is not circular: Petrow–Young arrive at a long fourth moment that can be bounded optimally via approximate functional equations and the spectral large sieve, whereas our short fourth moment can never be bounded optimally via the spectral large sieve, and it never implies subconvexity. We need to initially apply Motohashi's formula and then use the large sieve afterward.

Balancing in Corollary 1.3 leads one to a variant of Iwaniec's [35] short interval fourth moment bound.

**Corollary 1.4.** *Let $q \ll T^{1/2-\epsilon}$. For any primitive Dirichlet character $\chi$ modulo a prime $q$, we have*

$$\int_T^{T+(qT)^{2/3}} \left| L\left(\frac{1}{2}+it,\chi\right) \right|^4 dt \ll_\epsilon (qT)^{2/3+\epsilon}. \tag{1.9}$$

Jutila [42, Theorem 3] obtained the fourth moment bound

$$\sum_{\chi(\mathrm{mod}\,D)} \int_T^{T+T^{2/3}} \left| L\left(\frac{1}{2}+it,\chi\right) \right|^4 dt \ll_\epsilon D^{1+\epsilon} T^{2/3+\epsilon}. \tag{1.10}$$

This implies a bound for individual characters of the form $q^{1+\epsilon} T^{2/3+\epsilon}$. In the range $q \ll T^{1/2}$, our result improves upon this bound for individual characters to $(qT)^{2/3+\epsilon}$. There exists a fundamental obstruction that the parameter $H$ in Corollary 1.3 cannot exceed $T$, and the optimisation of $H$ in the bound in equation (1.8) in turn gives some restriction on $q$. Indeed, there is a weaker version where one has no restriction on the conductor, which can be deduced from the choice $H = T^{2/3}$ in Corollary 1.3. It therefore follows that

$$\int_T^{T+T^{2/3}} \left| L\left(\frac{1}{2}+it,\chi\right) \right|^4 dt \ll_\epsilon q^{1+\epsilon} T^{2/3+\epsilon}.$$

This is similar to equation (1.10) but eventually yields weaker subconvexity bounds for Dirichlet $L$-functions.

Corollary 1.4 is equivalent to the claim that Dirichlet $L$-functions cannot sustain large values, namely

$$\#\left\{ t \in \left[T, T+(qT)^{2/3}\right] : \left| L\left(\frac{1}{2}+it,\chi\right) \right| \geqslant V \right\} \ll_\epsilon (qT)^{2/3+\epsilon} V^{-4}.$$

It follows that $V \ll (qT)^{1/6+\epsilon}$, which yields Weyl-strength subconvexity bounds for Dirichlet $L$-functions when $q \ll (1+|t|)^{1/2-\epsilon}$ (the exponent $1/6$ often reoccurs in modern incarnations of these problems):

$$L\left(\frac{1}{2}+it,\chi\right) \ll_\epsilon (q(1+|t|))^{1/6+\epsilon}.$$

This is covered in Petrow–Young [63, 64], who have shown Weyl-strength subconvexity bounds for twisted automorphic $L$-functions, which are hybrid in both the $t$- and the $q$-aspect. We also establish the following:

**Corollary 1.5.** *Assume that $(qT)^{1/2+\epsilon} \ll T_0 \ll (qT)^{2/3}$ with $q \leqslant T^{1/2-\epsilon}$ and $T \leqslant t_1 < \cdots < t_R \leqslant 2T$ with $t_{r+1} - t_r \geqslant T_0$. For any primitive Dirichlet character $\chi$ modulo a prime $q$, we then have that*

$$\sum_{r=1}^R \int_{t_r}^{t_r+T_0} \left| L\left(\frac{1}{2}+it,\chi\right) \right|^4 dt \ll_\epsilon \left( RT_0 + qT\sqrt{\frac{R}{T_0}} \right)(qT)^\epsilon.$$

This kind of fourth moment bound is occasionally seen in the antecedent literature; we mention work of Iwaniec [35] and Jutila–Motohashi [41] to name two. Since one wants to prove Corollary 1.5 as an application of Theorem 1.1, one exploits the method outlined in the Appendix in Jutila–Motohashi [41]. They were unable to show Motohashi's formula for the fourth moment of Dirichlet $L$-functions averaged over primitive Dirichlet characters, but for the second moment of Estermann zeta functions. The proof of Corollary 1.5 is nearly identical to that of Corollary 1.3; however, one must detect cancellations in the spectral sum, which was executed by Ivić–Motohashi [34]. We mention that Jutila–Motohashi [41] established with the same notation as in Corollary 1.5 that

$$\sum_{\chi \pmod{D}} \sum_{r=1}^{R} \int_{t_r}^{t_r+T_0} \left| L\left(\frac{1}{2}+it,\chi\right) \right|^4 dt \ll_\epsilon \left( DRT_0 + D^{2(1+\vartheta)/3}(RT)^{2/3} \right) T^\epsilon$$

with $0 \leqslant \vartheta \leqslant 7/64$ denoting an admissible exponent toward the Ramanujan–Petersson conjecture.

**Corollary 1.6.** *Let $q \ll T^{1/2-\epsilon}$. For any primitive Dirichlet character $\chi$ modulo a prime $q$, we have*

$$\int_0^T \left| L\left(\frac{1}{2}+it,\chi\right) \right|^{12} dt \ll_\epsilon (qT)^{2+\epsilon}. \tag{1.11}$$

This extends work of Heath-Brown [29], who obtained the twelfth moment bound

$$\int_0^T \left| \zeta\left(\frac{1}{2}+it\right) \right|^{12} dt \ll_\epsilon T^{2+\epsilon}. \tag{1.12}$$

One should point out that Heath-Brown deduced equation (1.12) in an entirely different manner via the combination of van der Corput's method with the Halász–Montgomery inequality. Our method is fundamentally limited in this regard due to the fact that an average over characters is not included; we cannot hope to prove subconvexity bounds for $L(1/2+it,\chi)$ in the $q$-aspect with $t$ fixed, since we are averaging over too small a family. The result that would be of most interest includes an average over Dirichlet characters modulo $q$, namely the conjectural upper bound

$$\sum_{\chi \pmod{q}} \int_0^T \left| L\left(\frac{1}{2}+it,\chi\right) \right|^{12} dt \ll_\epsilon (qT)^{2+\epsilon}.$$

This is a new proof of Weyl-strength subconvexity. If $q=1$, this is Heath-Brown's bound for the twelfth moment of the Riemann zeta function. If $T \ll 1$ and $q$ is a smooth squarefree modulus, this estimate is due to Nunes [60], although his method does not work for arbitrary $q$. Meurman [48] proved the weaker bound $O_\epsilon(q^{3+\epsilon}T^{2+\epsilon})$, which implies Weyl-strength subconvexity in the $t$-aspect but convexity in the $q$-aspect. Jutila–Motohashi [41] established that

$$\sum_{q \leqslant Q} \sum_{\chi \pmod{q}} \int_0^T \left| L\left(\frac{1}{2}+it,\chi\right) \right|^{12} dt \ll_\epsilon Q^{3+\epsilon}T^{2+\epsilon}.$$

Milićević–White [51] have studied this problem when $T \ll 1$ and $q$ is growing in the depth aspect.

**Remark 1.7.** It would be interesting to see what happens if we average the fourth moment of individual Dirichlet $L$-functions over primitive characters modulo $q$. By brute force, we deduce

$$\sum_{\chi \pmod{q}}^* \frac{\mathcal{H}_\pm(\chi,\psi)}{\tau(\overline{\chi})} = \psi(\pm 1)J(\psi,\psi)$$

and this is of absolute value $\sqrt{q}$, since $\psi$ is not quadratic. The final bound looks like

$$\sum_{\chi \,(\text{mod } q)}^{*} \int_{T}^{T+H} \left| L\left(\frac{1}{2}+it, \chi\right) \right|^4 dt \ll_{\epsilon} (Hq)^{1+\epsilon} + \left(\frac{T\sqrt{q}}{\sqrt{H}}\right)^{1+\epsilon}.$$

Note that this is not completely inconceivable as it does not give any subconvexity bounds, but only the convexity bound due to the presence of the main term. However, there is no hope of getting $(qT)^{2+\epsilon}$ for the twelfth moment, since the bound for the moment in a short interval does not produce subconvexity. In order to circumvent this drawback, the key idea is the following. In the $T$-aspect, the trick is to divide the integral into shorter intervals and bound each one individually. With an additional sum over Dirichlet characters modulo $q$, we break this sum into a sum over short cosets. A similar idea is underlying work of Petrow–Young [63]. We reserve such pursuits for future work.

Finally, Topaçoğullari [67] has manifested an asymptotic formula for the fourth moment of individual Dirichlet $L$-functions. By the specification of a test function $g$ in Theorem 1.1 and a sequence of standard manipulations such as a spectral large sieve, one may arrive at an asymptotic formula of the same quality.

**Corollary 1.8** (Topaçoğullari [67, Theorem 1.1]). *Let $\epsilon > 0$, which is not necessarily the same at each occurrence. Let $\chi$ be any primitive Dirichlet character modulo $q$. Then we have for $T \geqslant 1$ that*

$$\int_{0}^{T} \left| L\left(\frac{1}{2}+it, \chi\right) \right|^4 dt = \int_{0}^{T} P_{\chi}(\log t)\,dt + O_{\epsilon}(q^{2-3\vartheta}T^{1/2+\vartheta+\epsilon} + qT^{2/3+\epsilon}), \tag{1.13}$$

*where $P_{\chi}$ is a polynomial of degree 4 whose coefficients depend only on $q$.*

We omit the proof of Corollary 1.8 since our method is quite analogous to that of Topaçoğullari [67]. Improving upon the error term in equation (1.13) in the $q$-aspect requires some additional manoeuvres.

### 1.4. Sketch of modus operandi

We here render a heuristic overview of the genesis of the automorphic reciprocity shown in Theorem 1.1. This is a high-level sketch geared toward experts. For the sake of argument, we may rely on approximate tools, although our proof is based upon more precise inspection of $L$-functions in the region of absolute convergence. We must also ignore all correction factors, degenerate terms, polar terms and the $t$-average.

An astute reader understands our stratagem from the shape of Motohashi's formula in equation (1.5). The overall ideas are essentially inspired by Motohashi's seminal work [55, §4.3–4.7] that enables one to prove fairly explicit spectral identities. We would like to study what happens if we replace the Riemann zeta function $\zeta(1/2+it)$ in his argument with a Dirichlet $L$-function $L(1/2+it, \chi)$, and we observe that the presence of the Dirichlet character substantially complicates our analysis. In a single phrase, we first open the four zeta values as Dirichlet series, apply Atkinson's dissection and then handle the shifted convolution sums via an application of the Voronoĭ summation formula, followed by the Kloosterman summation formula (Kuznetsov formula) attached to Atkin–Lehner cusps. An initial shape of the off-diagonal term looks like

$$\sum_{n,m \asymp \sqrt{q}} \overline{\chi}(n)\chi(n+m)\tau(n)\tau(n+m) = \sum_{a,b\,(\text{mod } q)} \overline{\chi}(a)\chi(a+b) \sum_{\substack{n,m \asymp \sqrt{q} \\ n \equiv a\,(\text{mod } q) \\ m \equiv b\,(\text{mod } q)}} \tau(n)\tau(n+m). \tag{1.14}$$

This sum has naturally arisen in the antecedent literature such as [64, 67, 69], and $\tau(m)$ must be replaced with the divisor function $\sigma_{\lambda}(m)$ in our proof due to the convergence issue. One sifts out the congruence conditions on the right-hand side of equation (1.14) via primitive additive characters. We can use the orthogonality of the Ramanujan sum

$$\delta_{n \equiv a \,(\mathrm{mod}\ q)} = \frac{1}{q} \sum_{c \mid q} r_c(n - a).$$

In order to spectrally expand equation (1.14), it is necessary to separate the two variables in $\tau(n + m)$. We thus make use of the approximate functional equation for the divisor function

$$\sigma_\lambda(m) = \sum_{(\ell, q) = 1} \frac{S(m, 0; \ell)}{\ell^{1 - \lambda}} \varpi_\lambda\left(\frac{\ell}{\sqrt{m}}\right) + m^\lambda \sum_{(\ell, q) = 1} \frac{S(m, 0; \ell)}{\ell^{1 + \lambda}} \varpi_{-\lambda}\left(\frac{\ell}{\sqrt{m}}\right), \qquad (1.15)$$

where in the notation of Young [69], we set

$$\varpi_\lambda(x) = \frac{1}{2\pi i} \int_{(a)} x^{-w} \zeta^q(1 - \lambda + w) \frac{G(w)}{w} dw.$$

The formula in equation (1.15) is a simple alternative to the $\delta$-symbol method of Duke–Friedlander–Iwaniec [22] and plays a rôle in eliminating the pole of the Riemann zeta function $\zeta(s)$ appearing in the original Ramanujan expansion. In our actual proof, we create a zero that like a *deus ex machina* kills the pole from the Riemann zeta function. One then applies the $\mathrm{GL}_2$ Voronoĭ summation formula to the sum over $n$. In other words, the functional equation of the Estermann zeta function is used. Letting $q\overline{q} \equiv 1$ (mod $\ell$), $\ell\overline{\ell} \equiv 1$ (mod $q$) and $\tilde{\tau}(\psi)$ be the normalised Gauß sum, we are led to sums of the product of Kloosterman sums[3]

$$\sum_{(\ell, q) = 1} \frac{S(m\overline{q}, \pm n\overline{q}; \ell) S(a\overline{\ell}, \pm n\overline{\ell}; q)}{\ell} \approx \sum_{\psi \,(\mathrm{mod}\ q)} \overline{\psi}(\pm anm) \tilde{\tau}(\psi)^2 \sum_{(\ell, q) = 1} \psi(\ell)^2 \frac{S(m\overline{q}, \pm n\overline{q}; \ell)}{\ell}. \quad (1.16)$$

Here the second Kloosterman sum $S(a\overline{\ell}, \pm n\overline{\ell}; q)$ was encoded via the orthogonality relation for Dirichlet characters modulo $q$, and the character $\overline{\psi}(m)$ was also added on the right-hand side for technical brevity. Be aware of the great similarity between equation (1.16) and Motohashi's conjecture written down in [54]. One decomposes the $\psi$-sum into the sum over primitive nonquadratic characters and others, followed by the application of the transformation formula due to Blomer–Milićević [15]. We are in a position to use the Kloosterman summation formula of level at most $q^2$ and $(\infty, \infty)$ cusp-pair. In this way, we obtain three automorphic $L$-functions with an explicit calculation of the resulting sums over $m$ and $n$, thereby deriving

$$\int_{-\infty}^{\infty} \left| L\left(\frac{1}{2} + it, \chi\right) \right|^4 g(t) dt \longleftrightarrow \sum_{\pm} \sum_{d \mid q} \frac{\mu(d)}{d} \sum_{\psi \,(\mathrm{mod}\ q)}^{\#} \mathcal{H}_\pm(\chi, \psi)$$

$$\times \sum_{f \in \mathcal{B}_\kappa^*(dq, \psi^2)} \epsilon_f^{(1 \mp 1)/2} \frac{L(1/2, f \otimes \overline{\psi})^3}{L(1, \mathrm{Ad}^2 f)} \check{g}(t_f), \qquad (1.17)$$

where the character sum $\mathcal{H}_\pm(\chi, \psi)$ was defined in Section 1.2 and the transform $\check{g}$ involves the hypergeometric function and depends on $\pm$. The contribution of the quadratic characters is similarly described. A benefit of using Motohashi's classical method is that one can achieve an explicit formulation of $\check{g}$, which would be alluring from an aesthetic point of view. If $\psi$ is a primitive character modulo $q$ and $f \in \mathcal{B}_\kappa^*(dq, \psi^2)$, then $f \otimes \overline{\psi}$ has trivial central character and conductor dividing $q^2$ (Theorem A.1).

---

[3]In this sketch, we use the symbol $\approx$ to mean that the left-hand side may roughly be written as an expression resembling the right-hand side with an acceptable error term.

There is some structural beauty in our reciprocity equation (1.17), since we were forced to decompose the $\psi$-sum in terms of whether $\psi$ is primitive nonquadratic or not. The condition that $f \otimes \overline{\psi}$ has the trivial central character is decisive as we may rely on the result of Guo [27], which guarantees that $L(1/2, f \otimes \overline{\psi}) \geqslant 0$. We can then evaluate the cubic moment in equation (1.17) using a standard positivity argument.

There is also a noteworthy plan to contemplate the fourth moment of individual Dirichlet $L$-functions. Fix a primitive Dirichlet character $\psi$ modulo $q$. We take $b = 1$ in [10, Theorem 1], multiply $\mathcal{T}_{a,b,q}(s, u, v)$ by $\overline{\psi}(a)$ and then sum over $a \pmod{q}$ with the application of the orthogonality relation. As an aside, this process necessitates a little modification since $a$ is supposed to be squarefree there. Nevertheless, this assumption was only to keep their formulæ cleaner, and hence their method works more generally when $a$ is not squarefree. It is believed that the substantial simplification of [10, (1.14)] will happen. We remark that the argument in [10] relies essentially on the isobaric sum $4 = 3 + 1$ and dualises 3 afterward. Its chief novelty is the use of the twisted multiplicativity of Kloosterman sums, enabling us to circumvent the manipulation of Kloosterman sums associated to various Atkin–Lehner cusps.

**Remark 1.9.** This work is relevant to the sixth moment of the Riemann zeta function (see [43]) whose calculations were initially contained in this article, but the author decided to remove this part since they are heuristics.

### 1.5. Organisation of the article

We devote Section 2 to compiling a preparatory toolbox in particular the evaluation of various multiplicative functions followed by the $GL_2$ Voronoï summation formula and the Kloosterman summation formula. We also equip the reader with an exhaustive exposition of Kloosterman sums at various cusps. In Section 3, we prove Theorem 1.1 via a shifted convolution problem. In Section 4, we establish Corollaries 1.3, 1.4, 1.5 and 1.6. The methods involve a combination of classical analytic number theory and automorphic forms. Moreover, the size of the conductor for twists of Maaß newforms is evaluated in Appendix A.

### 1.6. Basic notation and conventions

Throughout this article, the letter $\epsilon$ represents an arbitrarily small positive quantity, not necessarily the same at each occurrence. An implicit constant may depend on $\epsilon$, but this will often be suppressed from the notation. The Vinogradov symbol $A \ll B$ or the big O notation $A = O(B)$ signifies that $|A| \leqslant C|B|$ for some constant $C$. We use the notation $e(\zeta) := \exp(2\pi i \zeta)$ and $e_\alpha(\zeta) := e(\zeta/\alpha)$ with $\zeta \in \mathbb{C}$ and $\alpha \in \mathbb{R}$. We assume that a test function $g$ satisfies the following assumptions:

**Convention 1.10.** The function $g$ is real valued on $\mathbb{R}$, and there exists a large positive constant $A$ such that $g(t)$ is holomorphic and $\ll (1+|t|)^{-A}$ on a sufficiently wide horizontal strip $|\Im(t)| \leqslant A$. All implicit constants in Vinogradov symbols and big O notation may possibly depend on $A$ (where applicable).

## 2. Arithmetic and automorphic tools

We compile background materials that we shall need afterward to establish Theorem 1.1. For starters, we prepare elementary lemmata that are suitable for evaluating the character sums arising in this work. Our focus is on the simplification of sums involving multiplicative characters, which are akin to [64, §6.1]. Second, we introduce the Estermann zeta function and demonstrate its properties such as a functional equation. Finally, we present automorphic machinery as well as a preliminary exposition of Kloosterman sums. An exhaustive account of the theoretical background is found in [9, 21, 23] and references therein.

### 2.1. Manipulations of character sums

The orthogonality relation asserts that

$$\sum_{a(\mathrm{mod}\ q)} \chi(a) = \begin{cases} \varphi(q) & \text{if } \chi = \chi_0, \\ 0 & \text{otherwise,} \end{cases} \qquad \sum_{\chi(\mathrm{mod}\ q)} \chi(a) = \begin{cases} \varphi(q) & \text{if } a \equiv 1 (\mathrm{mod}\ q), \\ 0 & \text{otherwise.} \end{cases} \tag{2.1}$$

For any Dirichlet character $\chi$ modulo $q$, let

$$\tau(\chi, h) = \sum_{b(\mathrm{mod}\ q)} \chi(b) e\left(\frac{bh}{q}\right) \tag{2.2}$$

denote the Gauß sum associated to characters on residue classes modulo $q$. One writes $\tau(\chi) = \tau(\chi, 1)$ as usual. Multiplying equation (2.2) by $\overline{\chi}(a)$ and summing over $\chi$, we derive via orthogonality equation (2.1) that

$$e\left(\frac{ah}{q}\right) = \frac{1}{\varphi(q)} \sum_{\chi(\mathrm{mod}\ q)} \overline{\chi}(a) \tau(\chi, h) \quad \text{if} \quad (a, q) = 1. \tag{2.3}$$

This is the Fourier expansion of additive characters in terms of the multiplicative ones. One can evaluate Gauß sums for general Dirichlet characters.

**Lemma 2.1.** *Let $\chi$ be a nontrivial Dirichlet character modulo $q$ induced by the primitive character $\chi^*$ modulo $q^*$. For an integer $n \geqslant 1$, we have that*

$$\tau(\chi, n) = \tau(\chi^*) \sum_{d | (n, q/q^*)} d \overline{\chi^*}\left(\frac{n}{d}\right) \chi^*\left(\frac{q}{dq^*}\right) \mu\left(\frac{q}{dq^*}\right).$$

*Proof.* See [38, Lemma 3.2], which is corrected in the list of errata on Kowalski's website. □

We indicate by $r_c(n)$ and $S(m, n; c)$ the Ramanujan sum and Kloosterman sum, respectively, as follows:

$$r_c(n) := \sum_{a(\mathrm{mod}\ c)}^{*} e\left(\frac{an}{c}\right) = \sum_{d|(n,c)} d\mu\left(\frac{c}{d}\right), \qquad S(m, n; c) := \sum_{d(\mathrm{mod}\ c)}^{*} e\left(\frac{md + n\overline{d}}{c}\right), \tag{2.4}$$

where the asterisk means that the summation is restricted to a reduced system of residues. We have the Weil bound

$$|S(m, n; c)| \leqslant (m, n, c)^{1/2} c^{1/2} \tau(c). \tag{2.5}$$

This gives the best possible bound for individual Kloosterman sums, whereas we are apt to make use of the Kuznetsov formula (Theorem 2.11) to obtain additional savings from the sum over the moduli $c$. Indeed, various arithmetic problems can be transformed into bounding sums of Kloosterman sums. The twisted multiplicativity for Kloosterman sums is occasionally exploited. We also use the Jacobi sum

$$J(\chi, \psi) = \sum_{a(\mathrm{mod}\ q)} \chi(a) \psi(1 - a).$$

In particular, when $\chi$ and $\psi$ are of the same modulus and $\chi\psi$ is primitive, the relation between the Gauß sum and Jacobi sum is illustrated as (see [38, (3.18)])

$$J(\chi, \psi) = \frac{\tau(\chi)\tau(\psi)}{\tau(\chi\psi)}. \tag{2.6}$$

In general, we can establish the following lemma:

**Lemma 2.2.** *Assume that $q$ is prime. Suppose that $\psi_1, \psi_2$ are primitive characters modulo $q$ satisfying $\psi_1 \neq \overline{\psi_2}$, and let $a, b, c, d \in \mathbb{Z}$ with $(a, c, q) = 1$. Then*

$$\sum_{t \, (\mathrm{mod}\ q)} \psi_1(at + b)\psi_2(ct + d) = \overline{\psi_1}(c)\overline{\psi_2}(a)\psi_1\psi_2(ad - bc)\frac{\tau(\psi_1)\tau(\overline{\psi_1\psi_2})}{\tau(\overline{\psi_2})}. \qquad (2.7)$$

*Moreover, when $\psi_1 = \overline{\psi_2} = \psi$, we have that*

$$\sum_{t \, (\mathrm{mod}\ q)} \psi(at + b)\overline{\psi}(ct + d) = \psi(a)\overline{\psi}(c)r_q(ad - bc). \qquad (2.8)$$

*Proof.* The sum vanishes unless $(a, q) = (c, q) = 1$ as in [64, Page 455]. We assume $(a, q) = (c, q) = 1$ in what follows. To prove the first assertion, we initially expand the character $\psi_1$ into exponentials:

$$\psi_1(at + b) = \frac{1}{\tau(\overline{\psi_1})} \sum_{x \, (\mathrm{mod}\ q)} \overline{\psi_1}(x)e_q((at + b)x).$$

We would like to calculate the $t$-sum to reduce the problem to the manipulation of the $x$-sum. Then

$$\sum_{t \, (\mathrm{mod}\ q)} \psi_2(ct + d)e_q(atx) = \sum_{t \, (\mathrm{mod}\ q)} \psi_2(t)e_q(a\overline{c}(t - d)x) = \overline{\psi_2}(ax)\psi_2(c)e_q(-a\overline{c}dx)\tau(\psi_2).$$

On the other hand, one has

$$\sum_{x \, (\mathrm{mod}\ q)} \overline{\psi_1\psi_2}(x)e_q(x(b - a\overline{c}d)) = \psi_1\psi_2(ad - bc)\overline{\psi_1\psi_2}(-c)\tau(\overline{\psi_1\psi_2}).$$

Gathering the above identities together, we arrive at

$$\sum_{t \, (\mathrm{mod}\ q)} \psi_1(at + b)\psi_2(ct + d) = \overline{\psi_1}(c)\overline{\psi_2}(a)\psi_1\psi_2(bc - ad)\frac{\tau(\psi_2)\tau(\overline{\psi_1\psi_2})}{\tau(\overline{\psi_1})}.$$

Applying the relations $\tau(\overline{\psi_1}) = \psi_1(-1)\tau(\psi_1)^{-1}q$ and $\tau(\psi_2) = \psi_2(-1)\tau(\overline{\psi_2})^{-1}q$, the desired expression follows. The second claim is the same as in [64, Lemma 6.3]. □

As mentioned in the introduction, we will encounter sums of the product of Kloosterman sums. One should resolve $\overline{\ell}$ inside the argument of the Kloosterman sum, and the following lemma is helpful when we expand $S(a\overline{\ell}, \pm n\overline{\ell}; q)$ into multiplicative characters.

**Lemma 2.3.** *Assume that $q \geqslant 1$ and $a, b \in \mathbb{Z}$. We write $a = a_0a'$ and $b = b_0b'$, where $a_0b_0 \mid q^\infty$ and $(a'b', q) = 1$. Then*

$$S(a, b; q) = \frac{1}{\varphi(q)} \sum_{\psi \, (\mathrm{mod}\ q)} \tau(\psi, a_0)\tau(\psi, b_0)\overline{\psi}(a'b'), \qquad (2.9)$$

*where the sum runs over all Dirichlet characters modulo $q$ and $\varphi(q)$ is Euler's totient function.*

*Proof.* We exploit equation (2.3) so that the Kloosterman sum $S(a, b; q)$ equals

$$\sideset{}{^*}\sum_{d(\mathrm{mod}\ q)} e\left(\frac{a_0 d + b_0 a' b' \overline{d}}{q}\right) = \frac{1}{\varphi(q)} \sideset{}{^*}\sum_{d(\mathrm{mod}\ q)} e\left(\frac{a_0 d}{q}\right) \sum_{\psi(\mathrm{mod}\ q)} \overline{\psi}(a' b' \overline{d}) \tau(\psi, b_0)$$

$$= \frac{1}{\varphi(q)} \sum_{\psi(\mathrm{mod}\ q)} \tau(\psi, a_0) \tau(\psi, b_0) \overline{\psi}(a' b').$$

This finishes the proof of the lemma. □

Direct corollaries of Lemma 2.3 include

**Corollary 2.4.** *Suppose that $q \geqslant 1$ and $a, b \in \mathbb{Z}$. For $(ab, q) = 1$, we then have*

$$S(a, b; q) = \frac{1}{\varphi(q)} \sum_{\psi(\mathrm{mod}\ q)} \tau(\psi)^2 \overline{\psi}(ab).$$

This idea of separation of variables traces back to celebrated work of Blomer–Milićević in [15]. The number of the Gauß sums is the crux in Lemma 2.3 as the $n$th power determines the hyper-Kloosterman sum of $n$ variables. A variation of Lemma 2.3 was employed by Petrow–Young [64, Lemma 8.8], where they handled the hyper-Kloosterman sum $\mathrm{Kl}_3(x, y, z; q)$.

### 2.2. Double Estermann zeta function

We introduce the divisor function

$$\sigma_w(m) = \sum_{d | m} d^w.$$

As far as we know, the proof of the Hecke relation for $\sigma_w(m)$ has not appeared in the antecedent literature.

**Lemma 2.5.** *The divisor function satisfies the following multiplicativity relation:*

$$\sigma_w(mn) = \sum_{c | (m, n)} \mu(c) c^w \sigma_w\left(\frac{m}{c}\right) \sigma_w\left(\frac{n}{c}\right). \tag{2.10}$$

*Proof.* If $m = p_1^{k_1} p_2^{k_2} \cdots p_r^{k_r}$, $n = p_1^{\ell_1} p_2^{\ell_2} \cdots p_r^{\ell_r}$ and such a formula is shown for powers of primes, then

$$\sigma_w(mn) = \sigma_w(p_1^{k_1 + \ell_1}) \sigma_w(p_2^{k_2 + \ell_2}) \cdots \sigma_w(p_r^{k_r + \ell_r})$$

$$= \prod_{i=1}^{r} \sum_{j=0}^{\min(k_i, \ell_i)} \mu(p_i^j) p_i^{jw} \sigma_w(p_i^{k_i - j}) \sigma_w(p_i^{\ell_i - j}) \tag{2.11}$$

$$= \sum_{c | (m, n)} \mu(c) c^w \sigma_w\left(\frac{m}{c}\right) \sigma_w\left(\frac{n}{c}\right).$$

By multiplicativity, it suffices to prove the result for $m = p^k$, $n = p^\ell$. The formula is obvious if $\min(k, \ell) = 0$, so we assume both the variables are at least 1. Setting $X = p^w$, we find that the right-hand side of equation (2.11) reads

$$\text{RHS} = \sigma_w(p^k)\sigma_w(p^\ell) - p^w\sigma_w(p^{k-1})\sigma_w(p^{\ell-1})$$

$$= (1 + X + \cdots + X^k)(1 + X + \cdots + X^\ell) - X(1 + X + \cdots + X^{k-1})(1 + X + \cdots + X^{\ell-1})$$

$$= \frac{X^{k+\ell+1} - 1}{X - 1} = (1 + p^w + p^{2w} + \cdots + p^{(k+\ell)w}) = \sigma_w(p^{k+\ell}).$$

This establishes Lemma 2.5. □

For positive integers $h, \ell$ with $(h, \ell) = 1$ and $\Re(s) > 1$, we define the double Estermann zeta function as

$$D_2\left(s, \lambda; \frac{h}{\ell}\right) = \sum_{n=1}^{\infty} \sigma_\lambda(n) e\left(\frac{nh}{\ell}\right) n^{-s}. \tag{2.12}$$

Most of the analytic properties of $D_2(s, \lambda; h/\ell)$ follow from the identity

$$D_2\left(s, \lambda; \frac{h}{\ell}\right) = \ell^{\lambda-2s} \sum_{a,b \,(\text{mod } \ell)} e\left(\frac{abh}{\ell}\right) \zeta\left(s, \frac{a}{\ell}\right) \zeta\left(s - \lambda, \frac{b}{\ell}\right),$$

where for $\alpha \in \mathbb{R}$ and $\Re(s) > 1$,

$$\zeta(s, \alpha) := \sum_{n+\alpha>0} (n + \alpha)^{-s}$$

is the Hurwitz zeta function. It has meromorphic continuation to the entire complex plane $\mathbb{C}$ with a simple pole at $s = 1$ of residue 1 and satisfies the functional equation

$$\zeta(s, \alpha) = \sum_{\pm} G^{\mp}(1 - s)\zeta^{(\pm\alpha)}(1 - s), \tag{2.13}$$

where $G^{\pm}(s) = (2\pi)^{-s}\Gamma(s)e(\pm s/4)$ and $\zeta^{(\alpha)}(s)$ is a meromorphic continuation of $\sum_{n\geq 1} e(\alpha n)n^{-s}$. For $\alpha \in \mathbb{Q}$, the formula in equation (2.13) is a restatement of the Poisson summation in residue classes. Hence, as a function of the single variable $s$, the Estermann zeta function in equation (2.12) also has meromorphic continuation to all $s \in \mathbb{C}$ with two simple poles at $s = 1$ and $s = 1 + \lambda$ with respective residues $\ell^{\lambda-1}\zeta(1 - \lambda)$ and $\ell^{-\lambda-1}\zeta(1 + \lambda)$, provided $\lambda \neq 0$. In the case of $\lambda = 0$, there is a double pole at $s = 1$, and the Laurent expansion is given by

$$D_2\left(s, 0; \frac{h}{\ell}\right) = \frac{1/\ell}{(s-1)^2} + \frac{2(\gamma - \log \ell)/\ell}{(s-1)} + \cdots,$$

where $\gamma$ is the Euler–Mascheroni constant. The functional equation for $D_2(s, \lambda; h/\ell)$ reads as follows:

**Theorem 2.6** ([55, Lemma 3.7]). *The Estermann zeta function satisfies the functional equation*

$$D_2\left(s, \lambda; \frac{h}{\ell}\right) = 2(2\pi)^{2s-\lambda-2}\ell^{1+\lambda-2s}\Gamma(1-s)\Gamma(1+\lambda-s)$$

$$\times \left[ D_2\left(1-s, -\lambda; \frac{\overline{h}}{\ell}\right) \cos\left(\frac{\pi\lambda}{2}\right) - D_2\left(1-s, -\lambda; -\frac{\overline{h}}{\ell}\right) \cos\left(\pi\left(s - \frac{\lambda}{2}\right)\right) \right], \tag{2.14}$$

*where $\overline{h}$ is the multiplicative inverse of $h$ modulo $\ell$: that is, $h\overline{h} \equiv 1 \,(\text{mod } \ell)$.*

The formula in equation (2.14) was first proven by Hecke and Estermann in connection with an integral representation of the Hurwitz zeta function. We already find fragments of the Kloosterman sum in equation (2.14), which incline us to apply the Kloosterman summation formula. The functional equation of the Estermann zeta function $D_2$ is essentially equivalent to Voronoï summation. Indeed, the Estermann zeta function involves a divisor function instead of Hecke eigenvalues, and this corresponds to the standard Eisenstein series.

### 2.3. Kloosterman sums at singular cusps

Our presentation of cusps and scaling matrices is inspired by [37]. We restrict our attention to cusps with respect to the Hecke congruence subgroup

$$\Gamma = \Gamma_0(q) := \left\{ \begin{pmatrix} a & b \\ c & d \end{pmatrix} \in \mathrm{SL}_2(\mathbb{Z}) : c \equiv 0 (\mathrm{mod}\ q) \right\}.$$

Let $q_0 \mid q$. In the following, the letter $\psi$ denotes a Dirichlet character modulo $q_0$ with $\kappa = (1 - \psi(-1))/2$ such that $\psi(-1) = (-1)^\kappa$. Then we extend $\psi$ via the identification

$$\psi\left( \begin{pmatrix} a & b \\ c & d \end{pmatrix} \right) = \psi(d) = \overline{\psi}(a).$$

The group $\Gamma$ acts transitively on $\mathbb{P}^1(\mathbb{Q})$ by fractional linear transformations. An element $\mathfrak{a} \in \mathbb{P}^1(\mathbb{Q})$ is called a cusp. Two cusps $\mathfrak{a}$ and $\mathfrak{b}$ are termed equivalent under $\Gamma$ if there exists $\gamma \in \Gamma$ satisfying $\mathfrak{a} = \gamma \mathfrak{b}$. We write $q = rs$ with $(r, s) = 1$ and $q_0 \mid s$. Then we call a cusp of the form $\mathfrak{a} = 1/r$ an Atkin–Lehner cusp. The Atkin–Lehner cusps are equivalent to $\infty$ under the Atkin–Lehner operator as in [46, Definition 2.4]. Motohashi [57] singled out the scaling matrices corresponding to the Atkin–Lehner operators and this is crucial since the Fourier coefficients around the cusp $1/r$ are proportional to the Fourier coefficients around the cusp $\infty$ for an eigenform of the Atkin–Lehner operators.

Let $\Gamma_\mathfrak{a} = \{ \gamma \in \Gamma : \gamma \mathfrak{a} = \mathfrak{a} \}$ be the stabiliser of the cusp $\mathfrak{a}$ in $\Gamma$. A matrix $\sigma_\mathfrak{a} \in \mathrm{SL}_2(\mathbb{R})$, satisfying

$$\sigma_\mathfrak{a} \infty = \mathfrak{a} \quad \text{and} \quad \sigma_\mathfrak{a}^{-1} \Gamma_\mathfrak{a} \sigma_\mathfrak{a} = \Gamma_\infty = \left\{ \pm \begin{pmatrix} 1 & n \\ 0 & 1 \end{pmatrix} : n \in \mathbb{Z} \right\}$$

is called a scaling matrix for the cusp $\mathfrak{a}$. Since the scaling matrix $\sigma_\mathfrak{a}$ is not uniquely determined, the choice of $\sigma_\mathfrak{a}$ will be important in our subsequent discussions.

**Definition 2.7.** Let $\mathfrak{a}, \mathfrak{b}$ be cusps and $\sigma_\mathfrak{a}, \sigma_\mathfrak{b}$ be scaling matrices. The set

$$\mathcal{C}(\mathfrak{a}, \mathfrak{b}) = \left\{ c > 0 : \begin{pmatrix} * & * \\ c & * \end{pmatrix} \in \sigma_\mathfrak{a}^{-1} \Gamma \sigma_\mathfrak{b} \right\}$$

is called the set of allowed moduli.

For a cusp $\mathfrak{a}$ and a scaling matrix $\sigma_\mathfrak{a}$, let $u_\mathfrak{a}$ be such that $\sigma_\mathfrak{a}^{-1} u_\mathfrak{a} \sigma_\mathfrak{a} = \left( \begin{smallmatrix} 1 & 1 \\ & 1 \end{smallmatrix} \right)$. If $\psi$ satisfies $\psi(u_\mathfrak{a}) = 1$, then we say that $\mathfrak{a}$ is *singular* with respect to $\psi$. It behooves one to mention the following proposition:

**Proposition 2.8** ([46, Proposition 2.6]). *Assume that $q = rs$ with $(r, s) = 1$ with $q_0 \mid s$, where $q_0$ is the modulus of $\psi$. The two cusps $\infty$ and $1/r$ are then singular with respect to $\psi$. We choose a scaling matrix $\sigma_{1/r}$ associated to the Atkin–Lehner cusp $1/r$ to be an Atkin–Lehner operator, namely*

$$\sigma_{1/r} = \tau_r \nu_s \quad \text{with} \quad \tau_r = \begin{pmatrix} 1 & (s\overline{s} - 1)/r \\ r & s\overline{s} \end{pmatrix}, \quad \nu_s = \begin{pmatrix} \sqrt{s} & \\ & 1/\sqrt{s} \end{pmatrix}$$

*with $s\overline{s} \equiv 1 (\mathrm{mod}\ r)$. Then the set of allowed moduli is given by*

$$\mathcal{C}(\infty, 1/r) = \{ \gamma = c\sqrt{s} : c \equiv 0 (\mathrm{mod}\ r),\ (c, s) = 1 \}.$$

An important example is when $r = 1$ and $s = q$. In this case, one has

$$\mathcal{C}(\infty, 0) = \{c\sqrt{q} : c \geqslant 1, \ (c, q) = 1\}.$$

We now define Kloosterman sums with respect to a pair of cusps and general central character.

**Definition 2.9.** If $\mathfrak{a}, \mathfrak{b}$ are singular cusps for $\psi$ modulo $q_0$, then the Kloosterman sum associated to $\mathfrak{a}, \mathfrak{b}$ and $\psi$ with modulus $c$ is defined as

$$S_{\mathfrak{a}\mathfrak{b}}(m, n; c; \psi) = \sum_{\gamma = \left(\begin{smallmatrix} a & b \\ c & d \end{smallmatrix}\right) \in \Gamma_\infty \backslash \sigma_\mathfrak{a}^{-1} \Gamma \sigma_\mathfrak{b} / \Gamma_\infty} \psi(\operatorname{sgn}(c)) \overline{\psi(\sigma_\mathfrak{a} \gamma \sigma_\mathfrak{b}^{-1})} e\left(\frac{am + dn}{c}\right). \tag{2.15}$$

The occurrence of $-I \in \Gamma_\infty$ indicates that the lower-left entry $c$ is only defined up to $\pm$ sign, so that the factor $\psi(\operatorname{sgn}(c))$ accounts for this. Definition 2.9 is sensitively dependent upon the choice of $\sigma_\mathfrak{a}$ and $\sigma_\mathfrak{b}$. Moreover, if $|c| \notin \mathcal{C}(\mathfrak{a}, \mathfrak{b})$, then the sum appearing in equation (2.15) is empty, thus $S_{\mathfrak{a}\mathfrak{b}}(m, n; c; \psi) = 0$. If we stick to the case of $\mathfrak{a} = \infty$ and $\mathfrak{b} = 0$, we know the identity (see [46, (2.20)])

$$S_{\infty 0}(m, n; c\sqrt{q}; \psi) = \overline{\psi}(c) S(m, n\overline{q}; c) \tag{2.16}$$

with $(c, q) = 1$ and $q\overline{q} \equiv 1 \pmod{c}$. This differs from the one shown in [36, Page 58] by an additive character, which is due to a different choice of the scaling matrix. In the formula in equation (2.16), the presence of the character $\overline{\psi}(c)$ is a nice feature of the $(\infty, 0)$ cusp-pair as opposed to the $(\infty, \infty)$ cusp-pair:

$$S_{\infty \infty}(m, n; c; \psi) = S_\psi(m, n; c) = \sum_{d \,(\mathrm{mod}\ c)}^* \psi(d) e\left(\frac{md + n\overline{d}}{c}\right). \tag{2.17}$$

### 2.4. Normalisation and Hecke eigenvalues

In order to formulate our Kloosterman summation formula, we borrow the normalisation from [9, 23]. We refer the reader to [65] for an exposition in the case of general multiplier systems. Let $\Gamma \backslash \mathbb{H}$ be the modular surface, where $\Gamma = \Gamma_0(q)$ is the Hecke congruence subgroup and $\mathbb{H} = \{z \in \mathbb{C} : \Im(z) > 0\}$ is the upper half-plane. There exist various self-adjoint and pairwise commuting operators acting on the space $L^2(\Gamma \backslash \mathbb{H})$: the hyperbolic Laplacian

$$\Delta = -y^2\left(\frac{\partial^2}{\partial x^2} + \frac{\partial^2}{\partial y^2}\right) + i\kappa y \frac{\partial}{\partial x},$$

the Hecke operators (the non-archimedean counterparts of $\Delta$)

$$(T_n f)(z) = \frac{1}{\sqrt{n}} \sum_{ad=n} \psi(a) \sum_{b \,(\mathrm{mod}\ d)} f\left(\frac{az + b}{d}\right), \tag{2.18}$$

and the reflection operator $(T_{-1}f)(z) = f(-\overline{z})$, which flips positive and negative Fourier coefficients. As shall be explained later, we have the spectral decomposition of $L^2(\Gamma \backslash \mathbb{H})$ in terms of pure point spectrum, residual spectrum and continuous spectrum, namely

$$L^2(\Gamma \backslash \mathbb{H}) = L^2_{\mathrm{cusp}}(\Gamma \backslash \mathbb{H}) \oplus L^2_{\mathrm{res}}(\Gamma \backslash \mathbb{H}) \oplus L^2_{\mathrm{cont}}(\Gamma \backslash \mathbb{H}).$$

For an integer $k > 0$ with $k \equiv \kappa \pmod 2$, we choose a basis $\mathcal{B}_k(q, \psi)$ of holomorphic cusp forms. It is taken orthonormal with respect to the weight $k$ Petersson inner product

$$\langle h_1, h_2 \rangle = \int_{\Gamma \backslash \mathbb{H}} h_1(z) \overline{h_2(z)} y^k \frac{dx\,dy}{y^2}, \tag{2.19}$$

where $z = x + iy$. We let $\mathcal{B}_\kappa(q, \psi)$ denote a basis of the space of Maaß cusp forms. In particular, they are eigenfunctions on $\mathbb{H}$, are automorphic of weight $\kappa \in \{0, 1\}$, are square-integrable on a fundamental domain and vanish at all the cusps. Moreover, they are eigenfunctions of the $L^2$-extension of the hyperbolic Laplacian $\Delta$. For $f \in \mathcal{B}_\kappa(q, \psi)$, we write $\Delta f = s(1 - s)f$ with $s = 1/2 + it_f$ and $t_f \in \mathbb{R} \cup i[-1/2, 1/2]$. One may choose the basis $\mathcal{B}_\kappa(q, \psi)$ orthonormal with respect to the weight zero Petersson inner product introduced above. We define

$$\vartheta := \sup_{f \in \mathcal{B}_\kappa(q, \psi)} |\mathfrak{I} t_f|.$$

Then the Selberg eigenvalue conjecture asserts $\vartheta = 0$, whereas Selberg only established the upper bound $\vartheta \leqslant 1/4$. The current world record is $\vartheta \leqslant 7/64$ due to Kim–Sarnak [45] (see Blomer–Brumley [7] for the treatment in a more general scenario). The decomposition of the space of square-integrable, weight $\kappa$ automorphic forms on $\mathbb{H}$ with respect to the eigenspaces of the hyperbolic Laplacian involves the Eisenstein spectrum $\mathcal{E}(q, \psi)$, which is the orthogonal complement to the space of Maaß cusp forms. It is explicitly described in terms of the Eisenstein series $E_\mathfrak{a}(z, 1/2 + it)$, where $\mathfrak{a}$ runs over singular cusps and $t \in \mathbb{R}$. In this article, instead of using the classical Eisenstein series indexed by singular cusps as a basis of the continuous spectrum, we use another basis of Eisenstein series indexed by a set of parameters of the form

$$\{(\psi_1, \psi_2, f) : \psi_1 \psi_2 = \psi, \ f \in \mathcal{B}(\psi_1, \psi_2)\},$$

where $(\psi_1, \psi_2)$ ranges over the pairs of characters of modulus $q$ such that $\psi_1 \psi_2 = \psi$ and $\mathcal{B}(\psi_1, \psi_2)$ is a certain finite set dependent on $(\psi_1, \psi_2)$. We do not need to be more explicit here, and we refer the reader to [24] for an accurate definition of these parameters. The principal advantage of such a basis is that the Eisenstein series are eigenforms of the Hecke operators $T_n$ with $(n, q) = 1$: we have

$$T_n E(z, 1/2 + it, f) = \lambda_f(n, t) E(z, 1/2 + it, f)$$

with

$$\lambda_f(n, t) = \sum_{ab=n} \psi_1(a) a^{it} \psi_2(b) b^{-it}.$$

For convenience, we introduce the *ad hoc* notation

$$i(g, z) = cz + d, \qquad j(g, z) = (cz + d)|cz + d|^{-1}, \qquad g = \begin{pmatrix} * & * \\ c & d \end{pmatrix} \in \mathrm{SL}_2(\mathbb{R}).$$

We write the Fourier expansion of $f \in \mathcal{B}_k(q, \psi)$ around a singular cusp $\mathfrak{a}$ with a scaling matrix $\sigma_\mathfrak{a}$ as

$$i(\sigma_\mathfrak{a}, z)^{-k} f(\sigma_\mathfrak{a} z) = \left( \frac{(4\pi)^k}{\Gamma(k)} \right)^{1/2} \sum_{n=1}^\infty \rho_{f\mathfrak{a}}(n) n^{(k-1)/2} e(nz).$$

Let $f \in \mathcal{B}_\kappa(q, \psi)$ be an orthonormal basis of the space of Maaß cusp forms of weight $\kappa$ with respect to $\Gamma_0(q)$ and central character $\psi$. As is customary, we assume that each $f$ is either even or odd depending on whether $T_{-1}f = f$ or $T_{-1}f = -f$. If we denote the corresponding spectral parameter by $t_f$, we have

$$j(\sigma_{\mathfrak{a}}, z)^{-\kappa} f(\sigma_{\mathfrak{a}} z) = \sqrt{\cosh(\pi t_f)} \sum_{n \neq 0} \frac{\rho_{f\mathfrak{a}}(n)}{|n|^{1/2}} W_{\frac{n}{|n|} \frac{\kappa}{2}, it_f}(4\pi |n| y) e(nx)$$

with $W_{\alpha,\beta}$ the standard Whittaker function. For an Eisenstein series $E(z, 1/2 + it, f)$, we write

$$\begin{aligned} j(\sigma_{\mathfrak{a}}, z)^{-\kappa} E(\sigma_{\mathfrak{a}} z, 1/2 + it, f) = {} & c_{1,f,\mathfrak{a}}(t) y^{1/2+it} + c_{2,f,\mathfrak{a}}(t) y^{1/2-it} \\ & + \sqrt{\cosh(\pi t)} \sum_{n \neq 0} \frac{\rho_{f\mathfrak{a}}(n,t)}{|n|^{1/2}} W_{\frac{n}{|n|} \frac{\kappa}{2}, it}(4\pi |n| y) e(nx). \end{aligned}$$

Now let $\mathcal{B}_\kappa^*(q, \psi)$ be an orthonormal basis consisting of Hecke–Maaß newforms of level $q$ and central character $\psi \pmod{q_0}$ normalised so that $\lambda_f(1) = 1$, and we define similarly $\mathcal{B}_k^*(q, \psi)$. By Atkin–Lehner–Li theory [2, 3], we have the following direct sum decomposition as in [15]:

$$\mathcal{B}_\kappa(q, \psi) = \bigsqcup_{q_1 q_2 = q} \bigsqcup_{f \in \mathcal{B}_\kappa^*(q_1, \psi)} \mathcal{S}_{q_2}^{\square}(f), \tag{2.20}$$

where an element of the orthonormal basis $\mathcal{S}_{q_2}^{\square}(f) = \{f_{(d)}(z) : d \mid q_2\}$ is written as a linear combination of $f(cz)$ with $c \mid d$. A description of this linear combination is given by Iwaniec–Luo–Sarnak [39] for the squarefree level and principal nebentypus. Blomer–Milićević [14] subsequently orthonormalised the collection of Maaß forms $\{f(dz) : d \mid q_2\}$ for a newform $f$ and an integer $q_2 \geqslant 1$, which is not necessarily squarefree. Their scheme focuses on the principal nebentypus, but the generalisation to nontrivial central characters is feasible; see [31]. There are the Atkin–Lehner–Li theory for the continuous spectrum and a decomposition into spaces of oldforms analogous to equation (2.20) due to Young [70]. We make use of equation (2.20) to apply Hecke relations and the proportionality of Fourier coefficients to the Hecke eigenvalues.

In what follows, we make an abuse of notation $\rho_{f\infty}(n) = \rho_f(n)$ and $\rho_{f\infty}(n,t) = \rho_f(n,t)$. We consider the Fourier coefficients $\rho_{f\mathfrak{a}}(n)$ in more detail, and we stick to the case of $\mathfrak{a} \sim \infty$. Let $f \in \mathcal{B}_\kappa(q, \psi)$ be any Hecke eigenform, and let $\lambda_f(n)$ denote the corresponding eigenvalue for $T_n$, namely $T_n f = \lambda_f(n) f$. Then we often use the Hecke multiplicativity relations for $(mn, q) = 1$:

$$\lambda_f(mn) = \sum_{d \mid (m,n)} \mu(d) \psi(d) \lambda_f\left(\frac{m}{d}\right) \lambda_f\left(\frac{n}{d}\right), \qquad \lambda_f(m) \lambda_f(n) = \sum_{d \mid (m,n)} \psi(d) \lambda_f\left(\frac{mn}{d^2}\right). \tag{2.21}$$

We see some structural beauty of equation (2.21) in Section 3.5 along with the formula $\lambda_f(n) = \psi(n) \overline{\lambda_f(n)}$ for $(n, q) = 1$. Notice that the relation in equation (2.21) is valid for all $m, n \geqslant 1$ if $f$ is a newform. We invoke the bounds for the Hecke eigenvalues: if $f$ belongs to $\mathcal{B}_k(q, \psi)$ or is an Eisenstein series, there follows that

$$|\lambda_f(n)| \leqslant \tau(n) \ll_\epsilon n^\epsilon$$

for any $\epsilon > 0$. For $f \in \mathcal{B}_\kappa(q, \psi)$, the general upper bound is available in the form

$$|\lambda_f(n)| \leqslant \tau(n) n^\vartheta \ll_\epsilon n^{\vartheta+\epsilon}. \tag{2.22}$$

There is a connection between the Fourier coefficients and the Hecke eigenvalues

$$\rho_f(n) = \rho_f(1) \lambda_f(n), \qquad |\rho_f(1)|^2 = \frac{\varphi(q)}{2q^2} \frac{1}{L(1, \mathrm{Ad}^2 f)} \tag{2.23}$$

for $(n, q) = 1$. We refer the reader to [23, (6.14)] for the former and to [6, (2.10)] for the latter. Here the adjoint square $L$-function is defined as

$$L(s, \mathrm{Ad}^2 f) = \zeta(2s) \sum_{n=1}^{\infty} \frac{\lambda_f(n^2)}{n^s}.  \qquad (2.24)$$

The formula in equation (2.23) is valid for all $n \geqslant 1$ if $f$ is a newform.

**Remark 2.10.** If $\pi$ is an automorphic representation of $\mathrm{GL}_n$ with a contradredient representation $\widetilde{\pi}$, then it holds that $L(s, \widetilde{\pi}) = \overline{L(\overline{s}, \pi)}$. In particular, the Dirichlet series coefficients of $L(s, \pi)$ are real if and only if $\pi \cong \widetilde{\pi}$. If $n = 2$, we have $\widetilde{\pi} \cong \pi \otimes \omega^{-1}$, where $\omega$ denotes the central character of $\pi$. Hence the Hecke eigenvalues are not necessarily real in our context.

The aforementioned formulæ for Hecke operators acting on holomorphic and Maaß cusp forms apply to the Eisenstein series since $\rho_f(n, t)$ are proportional to the Hecke eigenvalues $\lambda_f(n, t)$ (compare [23, (7.13)]). We also derive the relation

$$\rho_f(-n, t) = \frac{\Gamma(s + \kappa/2)}{\Gamma(s - \kappa/2)} \rho_f(n, t).$$

In Section 3, we are mainly interested in primitive Dirichlet characters $\psi^2$ for primitive nonquadratic characters $\psi$ modulo a prime. In this case, one observes that $\kappa(\psi^2) = 0$, obtaining in particular that

$$|\rho_f(n, t)|^2 = \frac{|\lambda_f(n, t)|^2}{q |L(1 + 2it, \overline{\psi_1 \psi_2})|^2}.$$

We also need to define the twisted automorphic $L$-function

$$L(s, f \otimes \psi) := \sum_{n=1}^{\infty} \frac{\psi(n) \lambda_f(n)}{n^s}, \qquad \Re(s) > 1,  \qquad (2.25)$$

which has meromorphic continuation to the whole complex plane $\mathbb{C}$.

## 2.5. A version of the Kuznetsov formula

A spectral summation formula that we deploy here is an asymmetric trace formula relating sums of Kloosterman sums to Fourier coefficients of automorphic forms in the sense that the spectral side and the geometric side have a quite different shape. This technology is able to establish that there are considerable cancellations in sums of Kloosterman sums. Moreover, it provides a separation of variables $m$ and $n$ in the sum $\sum_c S(m, n; c)/c$, which is conducive to obtaining additional savings in summations over $m$ and $n$.

For $x > 0$, we introduce the three integral kernels

$$\mathcal{H}^+(x, t) := \frac{2\pi i t^{\kappa}}{\sinh(\pi t)} (J_{2it}(x) - (-1)^{\kappa} J_{-2it}(x)),$$

$$\mathcal{H}^-(x, t) := \frac{2\pi i^{1-\kappa}}{\sinh(\pi t)} (I_{2it}(x) - I_{-2it}(x)) = 8i^{-\kappa} \cosh(\pi t) K_{2it}(x),$$

$$\mathcal{H}^{\mathrm{hol}}(x, k) := 4i^k J_{k-1}(x), \quad k \in 2\mathbb{N},$$

where we have borrowed the normalisation from [67]. We decided to suppress $\psi$ from the notation, and we adopt this kind of abuse of notation unless otherwise specified. For $F \in \mathcal{C}_c^\infty(\mathbb{R}_+)$, we introduce

$$\mathscr{L}^\diamond F = \int_0^\infty \mathcal{H}^\diamond(\eta, \cdot) F(\eta) \frac{d\eta}{\eta}$$

for $\diamond \in \{+, -, \text{hol}\}$. With the whole notation set up, we formulate the Kloosterman summation formula, which in our case reads as follows:

**Theorem 2.11** ([67, Theorem 3.2]). *Assume* $F \in C^3(0, \infty)$ *satisfies* $x^j F^{(j)}(x) \ll \min(x, x^{-3/2})$ *for* $0 \leqslant j \leqslant 3$. *Let* $\mathfrak{a}, \mathfrak{b}$ *be singular cusps and* $\psi \,(\mathrm{mod}\, q_0)$ *be a Dirichlet character, and let* $m, n \geqslant 1$. *Then we have that*

$$\mathcal{O}_{\mathfrak{a}\mathfrak{b}}^q(m, \pm n; F; \psi) = \mathcal{A}_{\mathfrak{a}\mathfrak{b}}^{\mathrm{Maa}}(m, \pm n; \mathscr{L}^\pm F; \psi) + \mathcal{A}_{\mathfrak{a}\mathfrak{b}}^{\mathrm{Eis}}(m, \pm n; \mathscr{L}^\pm F; \psi) + \delta_{\pm=+} \mathcal{A}_{\mathfrak{a}\mathfrak{b}}^{\mathrm{hol}}(m, n; \mathscr{L}^{\,\mathrm{hol}} F; \psi).$$

*Here we set*

$$\mathcal{O}_{\mathfrak{a}\mathfrak{b}}^q(m, \pm n; F; \psi) \coloneqq \sum_{c \in \mathcal{C}(\mathfrak{a}, \mathfrak{b})} \frac{S_{\mathfrak{a}\mathfrak{b}}(m, \pm n; c; \psi)}{c} F\left(\frac{4\pi\sqrt{mn}}{c}\right), \tag{2.26}$$

$$\mathcal{A}_{\mathfrak{a}\mathfrak{b}}^{\mathrm{Maa}}(m, \pm n; F; \psi) \coloneqq \sum_{f \in \mathcal{B}_\kappa(q, \psi)} \overline{\rho_{f\mathfrak{a}}(m)} \rho_{f\mathfrak{b}}(\pm n) F(t_f), \tag{2.27}$$

$$\mathcal{A}_{\mathfrak{a}\mathfrak{b}}^{\mathrm{Eis}}(m, \pm n; F; \psi) \coloneqq \sum_{\psi_1\psi_2=\psi} \sum_{f \in \mathcal{B}(\psi_1, \psi_2)} \int_{-\infty}^\infty \overline{\rho_{f\mathfrak{a}}(m, t)} \rho_{f\mathfrak{b}}(n, t) F(t) \frac{dt}{4\pi}, \tag{2.28}$$

$$\mathcal{A}_{\mathfrak{a}\mathfrak{b}}^{\mathrm{hol}}(m, n; F; \psi) \coloneqq \sum_{\substack{k > \kappa \\ k \equiv \kappa(\mathrm{mod}\, 2)}} \sum_{f \in \mathcal{B}_k(q, \psi)} \overline{\rho_{f\mathfrak{a}}(m)} \rho_{f\mathfrak{b}}(n) F(k), \tag{2.29}$$

*where the sum over* $c$ *on the right-hand side of equation (2.26) runs over all positive real numbers for which the Kloosterman sum* $S_{\mathfrak{a}\mathfrak{b}}(m, \pm n; c; \psi)$ *is non-empty.*

The above formulation mimics the Kuznetsov formula that Motohashi [55] utilised to show a spectral reciprocity for the smoothed fourth moment of the Riemann zeta function. The name of the Kloosterman summation formula stems from the fact that it expresses certain sums of Kloosterman sums weighted by a function $F$ in terms of certain sums over automorphic forms weighted by transformed functions $\mathscr{L}^\diamond F$ with $\diamond \in \{+, -, \text{hol}\}$. We emphasise that there is no delta term in the Kloosterman summation formula.

## 3. Proof of Theorem 1.1

In this section, we prove Theorem 1.1 by applying the triad of Ramanujan–Voronoï–Kuznetsov to attack the shifted divisor sum in equation (1.14) for a primitive Dirichlet character $\chi$ modulo a prime $q$.

### 3.1. Dissection argument of Atkinson

Recall equation (1.3) and Convention 1.10. Let $\mathcal{R}_4^+$ (respectively, $\mathcal{R}_4^-$) be the subdomain of $\mathbb{C}^4$, where all four parameters have real parts larger than (respectively, less than) one. To wit, we define

$$\mathcal{R}_4^+ \coloneqq \{(s_1, s_2, s_3, s_4) \in \mathbb{C}^4 : \Re(s_i) > 1,\ 1 \leqslant i \leqslant 4\},$$
$$\mathcal{R}_4^- \coloneqq \{(s_1, s_2, s_3, s_4) \in \mathbb{C}^4 : \Re(s_i) < 1,\ 1 \leqslant i \leqslant 4\}.$$

For convenience, we write $\mathcal{Z}_2 := \mathcal{Z}_2(\mathbf{s}; g; \chi) = \mathcal{Z}_2(s_1, s_2, s_3, s_4; g; \chi)$ for a vector $\mathbf{s} = (s_1, s_2, s_3, s_4)$. In the domain $\mathcal{R}_4^+$, we open up the Dirichlet series to recast $\mathcal{Z}_2$ as

$$\int_{-\infty}^{\infty} \sum_{n_1=1}^{\infty} \sum_{n_2=1}^{\infty} \sum_{n_3=1}^{\infty} \sum_{n_4=1}^{\infty} \frac{\chi(n_1 n_2)\overline{\chi}(n_3 n_4)}{n_1^{s_1} n_2^{s_2} n_3^{s_3} n_4^{s_4}} \left(\frac{n_1 n_2}{n_3 n_4}\right)^{-it} g(t)dt$$
$$= \sum_{m=1}^{\infty} \sum_{n=1}^{\infty} \frac{\chi(m)\overline{\chi}(n)\sigma_{s_1-s_2}(m)\sigma_{s_3-s_4}(n)}{m^{s_1} n^{s_3}} \hat{g}\left(\frac{1}{2\pi}\log\frac{m}{n}\right), \quad (3.1)$$

where $\hat{g}$ is the Fourier transform[4]

$$\hat{g}(\xi) = \int_{-\infty}^{\infty} g(t)e(-\xi t)dt.$$

Shifting the contour in equation (1.3) yields that $\mathcal{Z}_2$ is meromorphic over the domain

$$\mathcal{B}_4 = \{(s_1, s_2, s_3, s_4) \in \mathbb{C}^4 : |s_1|, |s_2|, |s_3|, |s_4| < B\},$$

where $B = cA$ and $c$ is a small positive constant so that $B$ is sufficiently large. Although it is possible to make $B$ tend to infinity, we deal with the regime $\mathcal{B}_4$ for technical convenience. Since we are assuming that $\chi$ is primitive, the fourth moment $\mathcal{Z}_2$ is holomorphic in the vicinity of the point $\mathbf{s} = (1/2, 1/2, 1/2, 1/2)$. As an initial manipulation, we apply Atkinson's dissection, decomposing the double sum in equation (3.1) into three parts. It therefore follows that

$$\mathcal{Z}_2 = \mathcal{D}_4 + \mathcal{OD}_4^{\dagger} + \mathcal{OD}_4^{\ddagger}, \quad (3.2)$$

where $\mathcal{D}_4$ is the diagonal contribution corresponding to the terms with $m = n$, whereas $\mathcal{OD}_4^{\dagger}$ (respectively, $\mathcal{OD}_4^{\ddagger}$) is the off-diagonal contribution corresponding to the terms with $m > n$ (respectively, $m < n$). Note that the case of $m < n$ is symmetric to the case of $m > n$ as follows:

$$\mathcal{OD}_4^{\ddagger}(s_1, s_2, s_3, s_4; g; \chi) = \overline{\mathcal{OD}_4^{\dagger}(\overline{s_3}, \overline{s_4}, \overline{s_1}, \overline{s_2}; g; \chi)}. \quad (3.3)$$

Breaking the summation in this manner leads one to certain sums of Kloosterman sums and ultimately to automorphic forms. For an explanation of the reason this is not as surprising as it initially seems, we refer the reader to the discussion in [55, §4.2].

**Lemma 3.1.** *With the notation as above, we have that*

$$\mathcal{D}_4 = \hat{g}(0)\frac{\zeta^q(s_1+s_3)\zeta^q(s_1+s_4)\zeta^q(s_2+s_3)\zeta^q(s_2+s_4)}{\zeta^q(s_1+s_2+s_3+s_4)},$$

*where, for any L-function, the notation $L^q$ signifies the removal of the Euler factor at primes dividing $q$.*

*Proof.* It suffices to calculate

$$\mathcal{D}_4 = \hat{g}(0) \sum_{(m,q)=1} \frac{\sigma_{s_1-s_2}(m)\sigma_{s_3-s_4}(m)}{m^{s_1+s_3}},$$

where the Dirichlet series coefficient is a multiplicative function. Since the sum is over positive integers coprime to $q$, one can apply Ramanujan's identity (see [28, §17.8, Theorem 305])

$$\sum_{m=1}^{\infty} \frac{\sigma_{\alpha}(m)\sigma_{\beta}(m)}{m^s} = \frac{\zeta(s)\zeta(s-\alpha)\zeta(s-\beta)\zeta(s-\alpha-\beta)}{\zeta(2s-\alpha-\beta)},$$

omitting the Euler factors dividing $q$. This formula is valid as long as $\Re(s_i) > 1/2$. □

---

[4]A slightly different definition of the Fourier transform was used in Motohashi's monograph [55].

Note that $\hat{g}(0) = \int_{-\infty}^{\infty} g(t)dt$ is the mass of the weight function $\hat{g}_T(0) = T\hat{g}(0)$, where $g_T(x) = g(x/T)$ with $T > 0$. Hence, if $g$ is compactly supported on the interval $[1, 2]$, then $\mathrm{supp}(g_T) \subseteq [T, 2T]$. Since $\mathcal{D}_4$ has a pole of order 4 at $\mathsf{s} = (1/2, 1/2, 1/2, 1/2)$, we should have cancellation with a similar term coming from the off-diagonal terms. This is a common feature in the study of moment problems. The continuation of the off-diagonal contribution necessitates the full machinery of spectral theory of automorphic forms associated to Hecke congruence subgroups. We handle the second term in equation (3.2) due to the symmetry in equation (3.3), and the contribution from the third term will be incorporated at the end.

### 3.2. Evaluation of off-diagonal terms

For notational convenience, we call

$$G(y, s) = (1 + y)^{-s}\hat{g}\left(\frac{1}{2\pi}\log(1 + y)\right) = \int_{(1+\delta)} \mathring{g}(\tau, s)y^{-\tau}\frac{d\tau}{2\pi i} \qquad (3.4)$$

for suitable $\delta > 0$, where $\mathring{g}$ is the Mellin transform

$$\mathring{g}(\tau, s) = \int_0^{\infty} y^{\tau-1}(1 + y)^{-s}\hat{g}\left(\frac{1}{2\pi}\log(1 + y)\right)dy = \Gamma(\tau)\int_{-\infty}^{\infty}\frac{\Gamma(s - \tau + it)}{\Gamma(s + it)}g(t)dt, \qquad (3.5)$$

provided $\Re(s) > \Re(\tau) > 0$. The following lemma is useful:

**Lemma 3.2.** *As a function of two complex variables, $\mathring{g}(\tau, s)/\Gamma(\tau)$ is entire in $\tau$ and $s$. In addition, $\mathring{g}(\tau, s)$ is of rapid decay in $\tau$ as long as $s$ and $\Re(\tau)$ are bounded, say*

$$\mathring{g}(\tau, s) \ll (1 + |\tau|)^{-A}.$$

*Proof.* Shifting the contour $\Im(t) = 0$ in equation (3.5) downward appropriately, we obtain the first assertion. The second claim is a consequence of an upward shift. See [55, Lemma 4.1] for an analogous statement. $\square$

We now embark on the computation of $\mathcal{OD}_4^{\dagger}$ to resolve the shifted convolution problem. Pulling the Dirichlet characters out of the summations, one sees that

$$\mathcal{OD}_4^{\dagger} = \sum_{a,b(\mathrm{mod}\ q)} \overline{\chi}(a)\chi(a + b) \sum_{\substack{n \equiv a(\mathrm{mod}\ q) \\ m \equiv b(\mathrm{mod}\ q)}} \frac{\sigma_{s_3-s_4}(n)\sigma_{s_1-s_2}(n + m)}{n^{s_1+s_3}}G\left(\frac{m}{n}, s_1\right). \qquad (3.6)$$

We invoke the Ramanujan expansion (also known as the Ramanujan formula) for the divisor function, which is a precise formulation of equation (1.15). For any integers $n, q$ with $(n, q) = 1$ and $\Re(\xi) < 0$, we have that

$$\sigma_{\xi}(n) = \zeta^q(1 - \xi)\sum_{(\ell,q)=1}\ell^{\xi-1}r_{\ell}(n). \qquad (3.7)$$

This follows from the formula $r_{\ell}(n) = \sum_{d|(\ell,n)}d\mu(\ell/d)$ and a reversal of summations. Hence, we appeal to the fact that the divisor function emerges in the Fourier expansion of the Eisenstein series for $\mathrm{SL}_2(\mathbb{Z})$. On account of the presence of $\chi(a + b)$, we can assume $(n + m, q) = 1$. Upon applying equation (3.7), the right-hand side of equation (3.6) equals

$$\zeta^q(1 - s_1 + s_2)\sum_{a,b(\mathrm{mod}\ q)}\overline{\chi}(a)\chi(a + b)\sum_{(\ell,q)=1}\ell^{s_1-s_2-1}\sum_{\substack{n \equiv a(\mathrm{mod}\ q) \\ m \equiv b(\mathrm{mod}\ q)}}\frac{r_{\ell}(n + m)\sigma_{s_3-s_4}(n)}{n^{s_1+s_3}}G\left(\frac{m}{n}, s_1\right).$$

In order to simplify the innermost sums over $m$ and $n$, we detect the congruences modulo $q$ in an additive manner, obtaining

$$\mathcal{OD}_4^\dagger = \zeta^q(1 - s_1 + s_2)q^{-2} \sum_{a,b(\mathrm{mod}\ q)} \overline{\chi}(a)\chi(a+b) \sum_{(\ell,q)=1} \ell^{s_1-s_2-1}$$

$$\times \sum_{m=1}^{\infty} \sum_{n=1h(\mathrm{mod}\ \ell)}^{\infty} \sideset{}{^*}\sum \sum_{r|q\,c(\mathrm{mod}\ r)} \sideset{}{^*}\sum \sum_{d(\mathrm{mod}\ q)} e_r(c(n-a))e_q(d(m-b))e_\ell(h(n+m))\frac{\sigma_{s_3-s_4}(n)}{n^{s_1+s_3}}G\Big(\frac{m}{n}, s_1\Big).$$

The reason for the occurrence of the sum over $r$ is to make subsequent formulæ cleaner. If we execute the change of variables $(a, b) \mapsto (-a, -b)$, then the sum does not change according to whether $\chi$ is even or odd, and we use equation (3.4) to show that

$$\mathcal{OD}_4^\dagger = \zeta^q(1 - s_1 + s_2)q^{-2}$$

$$\times \sum_{a,b(\mathrm{mod}\ q)} \sum_{r|q} \sideset{}{^*}\sum_{c(\mathrm{mod}\ r)} \sum_{d(\mathrm{mod}\ q)} \overline{\chi}(a)\chi(a+b)e_r(ac)e_q(bd) \sum_{(\ell,q)=1} \ell^{s_1-s_2-1}$$

$$\times \sideset{}{^*}\sum_{h(\mathrm{mod}\ \ell)} \int_{(1+\delta)} \check{g}(\tau, s_1)\zeta^{(h/\ell+d/q)}(\tau)D_2\Big(s_1+s_3-\tau, s_3-s_4; \frac{h}{\ell}+\frac{c}{r}\Big)\frac{d\tau}{2\pi i}, \qquad (3.8)$$

where the right-hand side of equation (3.8) is absolutely convergent in a suitable domain such as

$$\mathcal{R}_{4,\delta} = \{(s_1, s_2, s_3, s_4) \in \mathcal{R}_4^+ : \Re(s_1+s_3) > 2 + 2\delta,\ \Re(s_1) + 1 < \Re(s_2),\ |s_3 - s_4| < \delta\}$$

for an arbitrary but fixed constant $\delta$. We then want to shift the contour to the right. In anticipation of the future application of the functional equation of the Estermann zeta function (Theorem 2.6), we define

$$\mathcal{E}_4 = \Big\{(s_1, s_2, s_3, s_4) \in \mathcal{B}_4 : \Re(s_1+s_3) < \frac{B}{3},\ \Re(s_1+s_4) < \frac{B}{3},\ \Re(s_1+s_2+s_3+s_4) > 3B\Big\}. \quad (3.9)$$

We ascertain that $\mathcal{R}_{4,\delta} \cap \mathcal{E}_4 \neq \emptyset$ if $\delta$ is small and confine the variables $(s_1, s_2, s_3, s_4)$ to be in the codomain $\mathcal{E}_4$, and then our shift of the contour in equation (3.8) to $\Re(\tau) = B$ results in

$$\sideset{}{^*}\sum_{h(\mathrm{mod}\ \ell)} \Big\{ \int_{(1+\delta)} - \int_{(B)} \Big\} \check{g}(\tau, s_1)\zeta^{(h/\ell+d/q)}(\tau)D_2\Big(s_1+s_3-\tau, s_3-s_4; \frac{h}{\ell}+\frac{c}{r}\Big)\frac{d\tau}{2\pi i}$$

$$= (\ell r)^{s_3-s_4-1}\zeta(1-s_3+s_4)\check{g}(s_1+s_3-1, s_1) \sum_{n=1}^{\infty} r_\ell(n)e_q(nd)n^{1-s_1-s_3}$$

$$+ (\ell r)^{s_4-s_3-1}\zeta(1+s_3-s_4)\check{g}(s_1+s_4-1, s_1) \sum_{n=1}^{\infty} r_\ell(n)e_q(nd)n^{1-s_1-s_4}. \quad (3.10)$$

Since the above polar terms do not contain the parameter $c$, one may express the sum over $c$ in equation (3.8) as $r_r(a)$, which leads one to a further simplification of them. At this stage, we notice that

$$\sum_{a,b(\mathrm{mod}\ q)} \overline{\chi}(a)\chi(a+b)e_q(ac+bd) = \tau(\chi, d)\overline{\tau(\chi, d-c)} = \overline{\chi}(d)\chi(d-c)q \qquad (3.11)$$

when $r = q$. This evaluation will not be used as the expanded sums are more amenable. Consequently, after some rearrangements, our formula is transformed into

$$\mathcal{OD}_4^\dagger = (\text{Polar Terms}) + q^{-2} \sum_{a,b,d(\text{mod } q)} \sum_{r|q} \zeta^r(1 - s_1 + s_2) \sum_{c(\text{mod } r)}^* \overline{\chi}(a)\chi(a+b)e_r(ac)e_q(bd)$$

$$\times \sum_{(\ell,r)=1} \ell^{s_1-s_2-1} \sum_{h(\text{mod } r)}^* \int_{(B)} \mathring{g}(\tau, s_1) \zeta^{(h/\ell+d/q)}(\tau) D_2\left(s_1 + s_3 - \tau, s_3 - s_4; \frac{h}{\ell} + \frac{c}{r}\right) \frac{d\tau}{2\pi i}$$

which is now in shape to make use of Poisson summation twice, namely Voronoĭ summation once. This is because we have $\Re(s_1 + s_3 - \tau) < 0$ and $\Re(s_1 + s_4 - \tau) < 0$; so we replace the Estermann zeta function $D_2(s_1 + s_3 - \tau, s_3 - s_4; h/\ell + c/r)$ with the absolutely convergent series exhibited by equation (2.14).

### 3.3. Utilisation of Voronoĭ summation

This subsection aims at applying the functional equation to the sum over $m$ and simplifying the resulting expression to a certain sum of Kloosterman sums of the form $S_{\mathfrak{ab}}(m, \pm n; c; \psi)$. To this end, the following elementary formula is necessary: if we set $v = hr + c\ell$ with $(v, \ell r) = 1$, then

$$\overline{v} = \overline{hr + c\ell} \equiv \overline{h}r\overline{r}^2 + \overline{c}\ell\overline{\ell}^2 (\text{mod } \ell r) \tag{3.12}$$

with $h\overline{h} \equiv 1$, $r\overline{r} \equiv 1(\text{mod } \ell)$ and $c\overline{c} \equiv 1$, $\ell\overline{\ell} \equiv 1(\text{mod } r)$. This holds under the assumption $(\ell, r) = 1$. The application of Theorem 2.6 hence leads one to

$$\mathcal{OD}_4^\dagger = (\text{Polar Terms}) + 2q^{-2} \sum_{a,b,d(\text{mod } q)} \sum_{r|q} \frac{\zeta^r(1 - s_1 + s_2)}{r} \sum_{c(\text{mod } r)}^* \overline{\chi}(a)\chi(a+b)e_r(ac)e_q(bd)$$

$$\times \sum_{(\ell,r)=1} \ell^{s_1-s_2-2} \sum_{h(\text{mod } \ell)}^* \int_{(B)} \left(\frac{2\pi}{\ell r}\right)^{2s_1+s_3+s_4-2\tau-2} \zeta^{(h/\ell+d/q)}(\tau)\Gamma(1 + \tau - s_1 - s_3)$$

$$\times \Gamma(1 + \tau - s_1 - s_4)\left\{D_2\left(1 + \tau - s_1 - s_3, s_4 - s_3; \frac{\overline{h}\overline{r}^2}{\ell} + \frac{\overline{c}\overline{\ell}^2}{r}\right)\cos\left(\frac{\pi(s_3 - s_4)}{2}\right)\right.$$

$$\left. - D_2\left(1 + \tau - s_1 - s_3, s_4 - s_3; -\frac{\overline{h}\overline{r}^2}{\ell} - \frac{\overline{c}\overline{\ell}^2}{r}\right)\cos\left(\pi\left(\frac{s_3 + s_4}{2} + s_1 - \tau\right)\right)\right\}\mathring{g}(\tau, s_1)\frac{d\tau}{2\pi i}.$$

We expand the integrand as the Dirichlet series once again so that the sums over $c$ and $h$ boil down to two Kloosterman sums, whence we are left with

$$\mathcal{OD}_4^\dagger = (\text{Polar Terms}) + 2(2\pi)^{s_1-s_2-1}q^{-2} \sum_{a,b,d(\text{mod } q)} \overline{\chi}(a)\chi(a+b)e_q(bd)$$

$$\times \sum_{r|q} r^{s_2-s_1}\zeta^r(1 - s_1 + s_2) \sum_{m=1}^\infty \sum_{n=1}^\infty m^{(1-s_1-s_2-s_3-s_4)/2}n^{(s_1-s_2+s_3-s_4-1)/2} \tag{3.13}$$

$$\times e_q(md)\sigma_{s_4-s_3}(n) \sum_\pm \sum_{(\ell,r)=1} \frac{1}{\ell}S(m, \pm n\overline{r}^2; \ell)S(a, \pm n\overline{\ell}^2; r)\Psi_{\mathsf{s}}^\pm\left(\frac{4\pi\sqrt{mn}}{\ell r}\right),$$

where

$$\Psi_{\mathsf{s}}^+(x) := \cos\left(\frac{\pi(s_3 - s_4)}{2}\right)\int_{(B)} \left(\frac{x}{2}\right)^{s_1+s_2+s_3+s_4-1-2\tau}\Gamma(1 + \tau - s_1 - s_3)\Gamma(1 + \tau - s_1 - s_4)\mathring{g}(\tau, s_1)\frac{d\tau}{2\pi i}, \tag{3.14}$$

$$\Psi_s^-(x) := -\int_{(B)} \left(\frac{x}{2}\right)^{s_1+s_2+s_3+s_4-1-2\tau} \cos\left(\pi\left(\frac{s_3+s_4}{2}+s_1-\tau\right)\right)$$
$$\times \Gamma(1+\tau-s_1-s_3)\Gamma(1+\tau-s_1-s_4)\mathring{g}(\tau,s_1)\frac{d\tau}{2\pi i}. \quad (3.15)$$

The expression in equation (3.13) affirmatively answers Motohashi's conjecture spelled out in his article [54] on the reciprocity formula for the second moment of Dedekind zeta functions. One ascertains that the integrand in equation (3.14) has rapid decay in $\tau$, whereas the definition in equation (3.15) requires the properties shown in Lemma 3.2. Convention 1.10 therefore applies to control $\Psi_s^-(x)$, which involves the opposite sign case of Kloosterman sums. This dominates the bulk of the evaluation of $\mathcal{OD}_4^\dagger$. We ensure this phenomenon in Proposition 3.5.

We then manage to simplify the product of the Kloosterman sums appearing in equation (3.13). It is desirable to collapse the second Kloosterman sum to perform the separation of variables. To this end, we write $n = n_0 n'$ with $n_0 \mid r^\infty$ and $(n', r) = 1$. We are in a position to use Lemma 2.3, getting

$$S(a, \pm n\overline{\ell}^2; r) = \frac{1}{\varphi(r)} \sum_{\psi \pmod r} \psi(\ell)^2 \overline{\psi}(\pm an')\tau(\psi)\tau(\psi, n_0), \quad (3.16)$$

where the condition $(a, r) = 1$ was used in view of the appearance of the character $\overline{\chi}(a)$ and $\tau(\psi)$ is the Gauß sum associated to $\psi$. Dropping the primes for notational simplicity, we observe that $\mathcal{OD}_4^\dagger$ equals the polar terms plus (note that the condition $(n, r) = 1$ is encoded by $\overline{\psi}(n)$)

$$2(2\pi)^{s_1-s_2-1}q^{-2} \sum_{a,b,d \pmod q} \overline{\chi}(a)\chi(a+b)e_q(bd) \sum_{r \mid q} \frac{r^{1-s_1+s_2}}{\varphi(r)} \zeta^r(1-s_1+s_2) \sum_{\pm} \sum_{\psi \pmod r}$$
$$\times \sum_{m=1}^\infty \sum_{n=1}^\infty \sum_{n_0 \mid r^\infty} \frac{\overline{\psi}(\pm an)\tau(\psi)\tau(\psi,n_0)e_q(md)\sigma_{s_4-s_3}(n_0 n)}{m^{(s_1+s_2+s_3+s_4-1)/2}(n_0 n)^{(1-s_1-s_2-s_3+s_4)/2}} \sum_{(\ell,r)=1} \psi(\ell)^2 \frac{S(m, \pm n_0 n\overline{r}^2; \ell)}{\ell r} \Psi_s^\pm\left(\frac{4\pi\sqrt{mn_0 n}}{\ell r}\right), \quad (3.17)$$

In order to adapt the Kloosterman summation formula in its suitable form, there are two roads to proceed. The first route is to represent $S(m, \pm n_0 n\overline{r}^2; \ell)$ by means of the Kloosterman sum associated to the $(\infty, 0)$ cusp-pair and apply Theorem 2.11 to the sum over $\ell$. The second is to extract the principal and quadratic character from the $\psi$-sum in equation (3.16) and then replace $S(m, \pm n_0 n\overline{r}^2; \ell)$ with the twisted Kloosterman sum; namely, we make use of the formula due to Blomer–Milićević [15, black(2.2)]:

$$\sum_{(\ell,r)=1} \psi(\ell)^2 \frac{S(m, \pm n_0 n\overline{r}^2; \ell)}{\ell r} \Psi_s^\pm\left(\frac{4\pi\sqrt{mn_0 n}}{\ell r}\right) = \frac{\psi(m)^2}{\tau(\psi^2)} \sum_{d \mid r} \mu(d) \sum_{dr \mid c} \frac{S_{\psi^2}(m, \pm n_0 n; c)}{c} \Psi_s^\pm\left(\frac{4\pi\sqrt{mn_0 n}}{c}\right), \quad (3.18)$$

where $S_\psi(m, n; c)$ is defined in equation (2.17) and $\psi$ denotes a primitive nonquadratic character modulo $r$ so that $\psi^2$ is primitive. This is where the assumption of $q$ being prime is used, since the square of a nonquadratic character is not necessarily primitive (see Remark 3.4). Nevertheless, one may remove this assumption via the use of an extended form of equation (3.18). Note that equation (3.18) is only available when $(m, r) = 1$. As shall be seen later, the Dirichlet character $\overline{\psi}(m)$ arises if we calculate the character sums in equation (3.17). We need to take the second route for technical simplicity. The off-diagonal term $\mathcal{OD}_4^\dagger$ is thus equal to the polar terms plus

$$2(2\pi)^{s_1-s_2-1}q^{-2}\sum_{a,b,c(\mathrm{mod}\,q)}\overline{\chi}(a)\chi(a+b)e_q(bc)\sum_{r|q}\frac{r^{1-s_1+s_2}}{\varphi(r)}\zeta^r(1-s_1+s_2)$$

$$\times\sum_{\pm}\left\{\sum_{\psi(\mathrm{mod}\,r)}^{\#}\frac{\overline{\psi}(\pm a)\tau(\psi)}{\tau(\psi^2)}\sum_{m=1}^{\infty}\sum_{n=1}^{\infty}\sum_{n_0|r^{\infty}}\frac{\psi(m)^2\overline{\psi}(n)\tau(\psi,n_0)e_q(mc)\sigma_{s_4-s_3}(n_0n)}{m^{(s_1+s_2+s_3+s_4-1)/2}(n_0n)^{(1-s_1+s_2-s_3+s_4)/2}}\right.$$

$$\left.\times\sum_{d|r}\mu(d)\sum_{dr|c}\frac{S_{\psi^2}(m,\pm n_0n;c)}{c}\Psi_{\mathsf{s}}^{\pm}\left(\frac{4\pi\sqrt{mn_0n}}{c}\right)+\mathcal{R}\right\},$$

where $\#$ on the $\psi$-sum means that the sum ranges over all primitive nonquadratic characters modulo $r$ and $\mathcal{R}$ serves as the remainder term

$$\mathcal{R}=\sum_{\psi(\mathrm{mod}\,r)}^{\flat}\sum_{m=1}^{\infty}\sum_{n=1}^{\infty}\sum_{n_0|r^{\infty}}\frac{\psi(\pm na)\tau(\psi)\tau(\psi,n_0)e_q(mc)\sigma_{s_4-s_3}(n_0n)}{m^{(s_1+s_2+s_3+s_4-1)/2}(n_0n)^{(1-s_1+s_2-s_3+s_4)/2}}$$

$$\times\sum_{(\ell,r)=1}\frac{S(m,\pm n_0n\overline{r}^2;\ell)}{\ell r}\Psi_{\mathsf{s}}^{\pm}\left(\frac{4\pi\sqrt{mn_0n}}{\ell r}\right).$$

Here $\flat$ on the $\psi$-sum means that the sum ranges over the principal and the quadratic character modulo $r$. Since it follows that $\tau(\psi,n_0)=\overline{\psi}(n_0)\tau(\psi)$ in the case where $\psi$ is a primitive nonquadratic character, we can take $n_0=1$ due to the condition $n_0\mid r^{\infty}$. In order to reduce the character sums over $a,b,c$, we start with simplifying the $c$-sum, obtaining

$$\sum_{a,b,c(\mathrm{mod}\,q)}\overline{\chi}(a)\overline{\psi}(a)\chi(a+b)e_q((b+m)c)=\psi(-1)q\sum_{a(\mathrm{mod}\,q)}\overline{\chi}(a)\overline{\psi}(a)\chi(a+m).$$

When $\chi\psi$ is primitive, one can set $\psi_1=\overline{\chi\psi}$ and $\psi_2=\chi$ in Lemma 2.2. To be more general, we deduce

$$\sum_{a(\mathrm{mod}\,q)}\overline{\chi}(a)\overline{\psi}(a)\chi(a+m)=\frac{1}{\tau(\overline{\chi})}\sum_{a,b(\mathrm{mod}\,q)}\overline{\chi}(ab)\overline{\psi}(a)e\left(\frac{(a+m)b}{q}\right)=\frac{\tau(\overline{\chi\psi})\tau(\psi,m)}{\tau(\overline{\chi})}.$$

Since the sum over primitive nonquadratic characters is empty when $r=1$, we have proven the following:

**Proposition 3.3.** *For any primitive Dirichlet character $\chi$ modulo a prime $q$, the function $\mathcal{OD}_4^{\dagger}$ can be meromorphically continued to the domain $\mathcal{E}_4$, and there we have*

$$\mathcal{OD}_4^{\dagger}=\mathcal{P}+\sum_{\pm}\{\mathcal{J}_{\pm}+\mathcal{E}_{\pm}\}.$$

*Here for the multiplicative function*

$$A_q(s_1,s_2,s_3,s_4):=\frac{1}{q}\sum_{c|q}\mu(c)c^{2-s_1-s_2-s_3-s_4}\sum_{d|q/c}\mu\left(\frac{q}{cd}\right)d^{2-s_1-s_3}\sigma_{s_1-s_2+s_3-s_4-1}(d),$$

*we set*

$$\mathcal{P}=A_q(s_1,s_2,s_3,s_4)\mathring{g}(s_1+s_3-1,s_1)\frac{\zeta^q(1-s_1+s_2)\zeta^q(1-s_3+s_4)\zeta(s_2+s_4)\zeta(s_1+s_3-1)}{\zeta^q(s_2+s_4-s_1-s_3+2)}$$

$$+A_q(s_1,s_2,s_4,s_3)\mathring{g}(s_1+s_4-1,s_1)\frac{\zeta^q(1-s_1+s_2)\zeta^q(1+s_3-s_4)\zeta(s_2+s_3)\zeta(s_1+s_4-1)}{\zeta^q(s_2+s_3-s_1-s_4+2)},$$

$$\mathcal{J}_{\pm} = \frac{2\zeta^q(1-s_1+s_2)}{\tau(\overline{\chi})\varphi(q)q}\left(\frac{2\pi}{q}\right)^{s_1-s_2-1}\sum_{\psi\,(\mathrm{mod}\,q)}^{\#}\psi(\mp 1)\tau(\psi)\tau(\overline{\chi\psi})J(\psi,\psi)\sum_{d|q}\mu(d)$$

$$\times\sum_{m=1}^{\infty}\sum_{n=1}^{\infty}\psi(m)\overline{\psi}(n)m^{(1-s_1-s_2-s_3-s_4)/2}n^{(s_1-s_2+s_3-s_4-1)/2}\sigma_{s_4-s_3}(n)\,\mathcal{O}_{\infty\infty0}^{dq}(m,\pm n;\Psi_s^{\pm};\psi^2),\quad(3.19)$$

$$\mathcal{E}_{\pm} = \frac{2(2\pi)^{s_1-s_2-1}}{\tau(\overline{\chi})q}\sum_{r|q}\frac{r^{1-s_1+s_2}}{\varphi(r)}\zeta^r(1-s_1+s_2)\sum_{\psi\,(\mathrm{mod}\,r)}^{\flat}\psi(\mp 1)\tau(\psi)\tau(\overline{\chi\psi})\sum_{m=1}^{\infty}\sum_{n=1}^{\infty}\sum_{n_0|r^\infty}\psi(n)$$

$$\times\tau(\psi,m)\tau(\psi,n_0)m^{(1-s_1-s_2-s_3-s_4)/2}(n_0n)^{(s_1-s_2+s_3-s_4-1)/2}\sigma_{s_4-s_3}(n_0n)\,\mathcal{O}_{\infty\infty0}^{r^2}(m,\pm n_0n;\Psi_s^{\pm};\psi_0),\quad(3.20)$$

*where $\mathcal{O}_{\mathrm{ab}}^q(m,\pm n;\Psi_s^{\pm};\psi)$ is defined in Theorem 2.11.*

**Remark 3.4.** The square of a primitive nonquadratic Dirichlet character of composite modulus having at least two prime factors is not primitive. For instance, one may choose a primitive quadratic character $\psi_1$ (respectively, nonquadratic character $\psi_2$) modulo 3 (respectively, modulo 5) and combine them to obtain a character $\psi_1\psi_2$ modulo 15 as follows:

| $n$ | 1 | 2 | 3 | 4 | 5 | 6 | 7 | 8 | 9 | 10 | 11 | 12 | 13 | 14 | 15 |
|---|---|---|---|---|---|---|---|---|---|---|---|---|---|---|---|
| $\psi_1(n)$ | 1 | $-1$ | 0 | 1 | $-1$ | 0 | 1 | $-1$ | 0 | 1 | $-1$ | 0 | 1 | $-1$ | 0 |
| $\psi_2(n)$ | 1 | $i$ | $-i$ | $-1$ | 0 | 1 | $i$ | $-i$ | $-1$ | 0 | 1 | $i$ | $-i$ | $-1$ | 0 |
| $\psi_1\psi_2(n)$ | 1 | $-i$ | 0 | $-1$ | 0 | 0 | $i$ | $i$ | 0 | 0 | $-1$ | 0 | $-i$ | 1 | 0 |

The multiplication of Dirichlet characters of different moduli gives a Dirichlet character modulo a least common multiple. We define $\psi=\psi_1\psi_2$, which is nonquadratic since it involves $i$ and $-i$ as values and is primitive since it cannot be induced from a Dirichlet character modulo 1, 3, or 5. However, $\psi^2=\psi_1^2\psi_2^2$ is induced from a quadratic character modulo 5. This is why we assume that $q$ is prime in this work.

*Proof.* For the first assertion, one can check that both series in equations (3.19) and (3.20) converge absolutely and uniformly in $\mathcal{E}_4$, proving the meromorphic continuation of $\mathcal{OD}_4^\dagger$ to $\mathbb{C}^4$. A remarkable point is that we only need the Weil bound in equation (2.5) for Kloosterman sums. For the second assertion, it suffices to compute the polar terms. We decompose $\mathcal{P}=\mathcal{P}_1+\mathcal{P}_2$, where $\mathcal{P}_1$ (respectively, $\mathcal{P}_2$) stems from the first term (respectively, second term) on the right-hand side of equation (3.10). First of all, we remark that the polar contribution $\mathcal{P}$ is expressed as

$$\zeta^q(1-s_1+s_2)q^{-2}\sum_{a,b\,(\mathrm{mod}\,q)}\sum_{r|q}\sum_{c\,(\mathrm{mod}\,r)}^{*}\sum_{d\,(\mathrm{mod}\,q)}\overline{\chi}(a)\chi(a+b)$$

$$\times e_r(ac)e_q(bd)\sum_{(\ell,q)=1}\ell^{s_1-s_2-1}\times(\text{RHS of equation (3.10)}).$$

In the following lines, we deal with $\mathcal{P}_1$. The sum over $\ell$ is reduced to

$$\sum_{(\ell,q)=1}\ell^{s_1-s_2+s_3-s_4-2}r_\ell(n)=\zeta^q(s_2+s_4-s_1-s_3+2)^{-1}\sigma_{s_1-s_2+s_3-s_4-1}(n).$$

We are in a position to manipulate the exponential sums. The sum over $d$ gives rise to

$$\sum_{d\,(\mathrm{mod}\,q)}e_q((b+n)d)=\begin{cases}q & \text{if }b\equiv -n\,(\mathrm{mod}\,q),\\0 & \text{otherwise.}\end{cases}$$

This renders

$$\sum_{a(\mathrm{mod}\ q)} \overline{\chi}(a)\chi(a-n)r_r(a) = \frac{\mu(r)}{\tau(\overline{\chi})}\sum_{a,b(\mathrm{mod}\ q)} \overline{\chi}(ab)e_q((a+n)b) = \mu(r)r_q(n).$$

In this manner, the sums over $a, b, c, d$ can be eliminated. In general, the sum over $r$ turns into

$$\sum_{r\mid q} \mu(r)r^{s_3-s_4-1} = \prod_{p\mid q}\left(1 - \frac{1}{p^{1-s_3+s_4}}\right),$$

which yields $\zeta^q(1-s_3+s_4)$ when combined with $\zeta(1-s_3+s_4)$. We hence obtain

$$\mathcal{P}_1 = \frac{\mathring{g}(s_1+s_3-1, s_1)}{q}\frac{\zeta^q(1-s_1+s_2)\zeta^q(1-s_3+s_4)}{\zeta^q(s_2+s_4-s_1-s_3+2)}\sum_{n=1}^{\infty}\frac{\sigma_{s_1-s_2+s_3-s_4-1}(n)r_q(n)}{n^{s_1+s_3-1}}.$$

Having mentioned this expression, we need to calculate the the sum over $n$ via Lemma 2.5, getting

$$\sum_{n=1}^{\infty}\frac{\sigma_w(n)}{n^s}\sum_{d\mid(n,q)}\mu\left(\frac{q}{d}\right)d = \sum_{d\mid q}\mu\left(\frac{q}{d}\right)d^{1-s}\sum_{n=1}^{\infty}\frac{1}{n^s}\sum_{c\mid(n,d)}\mu(c)\sigma_w\left(\frac{n}{c}\right)\sigma_w\left(\frac{d}{c}\right)c^w$$

$$= \sum_{d\mid q}\mu\left(\frac{q}{d}\right)d^{1-s}\sum_{c\mid d}\mu(c)c^{w-s}\sigma_w\left(\frac{d}{c}\right)\sum_{n=1}^{\infty}\frac{\sigma_w(n)}{n^s}$$

$$= \zeta(s)\zeta(s-w)\sum_{c\mid q}\mu(c)c^{1+w-2s}\sum_{d\mid q/c}\mu\left(\frac{q}{cd}\right)d^{1-s}\sigma_w(d),$$

where $s = s_1+s_3-1$ and $w = s_1-s_2+s_3-s_4-1$. Notice that the right-hand side can be written as a Dirichlet convolution of multiplicative functions. Recall that the Dirichlet convolution is defined as

$$(f*g)(m) = \sum_{d\mid m}f(d)g\left(\frac{m}{d}\right),$$

which is an associative and commutative binary operation on arithmetic functions. If $f$ and $g$ are assumed to be multiplicative, so is $f*g$. Setting $f(n) := \sigma_w(n)n^{1-s}$ and $g(n) := \mu(n)n^{1+w-2s}$, we have that

$$\sum_{n=1}^{\infty}\frac{\sigma_w(n)r_q(n)}{n^s} = (f*g*\mu)(q)\zeta(s)\zeta(s-w).$$

This finishes the proof of Proposition 3.3. □

### 3.4. Utilisation of the Kuznetsov formula

We now apply the automorphic machinery in Section 2 to our sums of Kloosterman sums along with a careful analysis of integrals involving various Bessel functions. We focus on the treatment of $\mathcal{O}^{dq}_{\infty\mathfrak{b}}(m, \pm n; \Psi^{\pm}_{\mathfrak{s}}; \psi)$ with $\mathfrak{b} = \infty, 0$ in an across-the-board manner and implement the substitutions at the end to calculate $\mathcal{J}_{\pm}$ and $\mathcal{E}_{\pm}$. The three-time differentiability and the decay condition at $x = 0$ in Theorem 2.11 are fulfilled. Moreover, the decay condition at infinity follows by shifting the contours $\Re(s) = B$ in equations (3.14) and (3.15) to the right and then applying Lemma 3.2. We forthwith derive the spectral decomposition

$$\mathcal{O}^{dr}_{\infty\mathfrak{b}}(m, \pm n; \Psi^{\pm}_{\mathfrak{s}}; \psi) = \mathcal{A}^{\mathrm{Maa\beta}}_{\infty\mathfrak{b}}(m, \pm n; \mathscr{L}^{\pm}\Psi^{\pm}_{\mathfrak{s}}; \psi) + \mathcal{A}^{\mathrm{Eis}}_{\infty\mathfrak{b}}(m, \pm n; \mathscr{L}^{\pm}\Psi^{\pm}_{\mathfrak{s}}; \psi)$$
$$+ \delta_{\pm=+}\mathcal{A}^{\mathrm{hol}}_{\infty\mathfrak{b}}(m, n; \mathscr{L}^{\mathrm{hol}}\Psi^{+}_{\mathfrak{s}}; \psi). \tag{3.21}$$

Complying with Motohashi's book [55], one observes how the case of the plus sign contributes. In order to reduce the integral transform $\mathscr{L}^+$, we consider, in the light of the definition in equation (3.14), the double integral

$$\int_0^\infty J_{2it}(\eta) \int_{(B)} \Gamma(1+\tau-s_1-s_3)\Gamma(1+\tau-s_1-s_4)\left(\frac{\eta}{2}\right)^{s_1+s_2+s_3+s_4-2-2\tau} \mathring{g}(\tau,s_1)\,d\tau\,d\eta, \qquad (3.22)$$

where $\mathring{g}$ is the Mellin transform in equation (3.5). This double integral is holomorphic in the domain

$$\left\{(s_1,s_2,s_3,s_4)\in\mathbb{C}^4: \begin{array}{c} \Re(s_1+s_3)<1+B,\ \Re(s_1+s_4)<1+B, \\ 1+2B<\Re(s_1+s_2+s_3+s_4)<3/2+2A \end{array}\right\}, \qquad (3.23)$$

which contains $\mathcal{E}_4$ defined in equation (3.9) and $A>B$ is the same as in Convention 1.10. To prove this fact, we divide equation (3.22) into two parts according as $0\leqslant\eta<1$ and $\eta\geqslant1$ and note that

$$J_{2it}(\eta)\ll\begin{cases}1 & \text{as } \eta\to0, \\ \eta^{-1/2} & \text{as } \eta\to\infty,\end{cases}$$

where we implicitly assumed that $t$ is real. The first part is clearly holomorphic in the domain of equation (3.23). To estimate the second part, we move the contour in the $\tau$-integral to $\Re(\tau)=A$. We would like to consider the subdomain of equation (3.23) where we have $1+2B<\Re(s_1+s_2+s_3+s_4)<3/2+2B$. There the integral in equation (3.22) is absolutely and uniformly convergent in view of the power series expansion of the $J$-Bessel function

$$J_\nu(y)=\sqrt{\frac{2}{\pi y}}\left\{\cos\left(y-\frac{\pi}{2}\left(\nu+\frac{1}{2}\right)\right)\sum_{m=0}^M(-1)^m\binom{\nu-1/2}{2m}\frac{\Gamma(\nu+1/2+2m)}{\Gamma(\nu+1/2)(2y)^{2m}}\right.$$
$$\left.-\sin\left(y-\frac{\pi}{2}\left(\nu+\frac{1}{2}\right)\right)\sum_{m=0}^M(-1)^m\binom{\nu-1/2}{2m+1}\frac{\Gamma(\nu+3/2+2m)}{\Gamma(\nu+1/2)(2y)^{2m+1}}\right\}+O(y^{-\Re(\nu)-3/2-2M}),$$

provided $\Re(\nu)>-2M-3/2$. Hence, the double integral in equation (3.22) equals

$$\int_{(B)}\frac{\Gamma\left(\frac{s_1+s_2+s_3+s_4-1}{2}-\tau+it\right)}{\Gamma\left(\frac{3-s_1-s_2-s_3-s_4}{2}+\tau+it\right)}\Gamma(1+\tau-s_1-s_3)\Gamma(1+\tau-s_1-s_4)\mathring{g}(\tau,s_1)\,d\tau, \qquad (3.24)$$

where we have used the formula

$$\int_0^\infty J_\nu(\eta)\left(\frac{\eta}{2}\right)^{-\mu}d\eta=\Gamma\left(\frac{1+\nu-\mu}{2}\right)\Gamma\left(\frac{1+\nu+\mu}{2}\right)^{-1}$$

valid for $1/2<\Re(\mu)<1+\Re(\nu)$. Since the integral in equation (3.24) is regular in equation (3.23), the analytic continuation implies that the double integral in equation (3.22) is equal to equation (3.24) throughout the domain of equation (3.23) or $\mathcal{E}_4$. At this stage, we note that the following identities hold:

$$\frac{\Gamma(a-\tau+it)}{\Gamma(1-a+\tau+it)}-\frac{\Gamma(a-\tau-it)}{\Gamma(1-a+\tau-it)}=\frac{2}{\pi i}\sinh(\pi t)\cos(\pi(a-\tau))\Gamma(a-\tau+it)\Gamma(a-\tau-it),$$
$$\frac{\Gamma(a-\tau+it)}{\Gamma(1-a+\tau+it)}+\frac{\Gamma(a-\tau-it)}{\Gamma(1-a+\tau-it)}=\frac{2}{\pi}\cosh(\pi t)\sin(\pi(a-\tau))\Gamma(a-\tau+it)\Gamma(a-\tau-it).$$

which follows from the functional equation $\Gamma(s)\Gamma(1-s) = \pi/\sin(\pi s)$. Substituting $a = (s_1 + s_2 + s_3 + s_4 - 1)/2$, we obtain after some rearrangement that[5]

$$\mathscr{L}^+\Psi_{\mathsf{s}}^+(t) = (it\coth(\pi t))^{\kappa}\cos\left(\frac{\pi(s_3-s_4)}{2}\right)\int_{(B)}\sin\left(\frac{\pi(s_1+s_2+s_3+s_4-\kappa-2\tau)}{2}\right)$$
$$\times\Gamma\left(\frac{s_1+s_2+s_3+s_4-1}{2}+it-\tau\right)\Gamma\left(\frac{s_1+s_2+s_3+s_4-1}{2}-it-\tau\right)$$
$$\times\Gamma(1+\tau-s_1-s_3)\Gamma(1+\tau-s_1-s_4)\mathring{g}(\tau,s_1)\frac{d\tau}{\pi i}.$$

In a similar fashion, the holomorphic contribution turns into

$$\mathscr{L}^{\mathrm{hol}}\Psi_{\mathsf{s}}^+(k) = i^{k-1}\cos\left(\frac{\pi(s_3-s_4)}{2}\right)\int_{(B)}\frac{\Gamma\left(\frac{k+s_1+s_2+s_3+s_4-2}{2}-\tau\right)}{\Gamma\left(\frac{k+2-s_1-s_2-s_3-s_4}{2}+\tau\right)}$$
$$\times\Gamma(1+\tau-s_1-s_3)\Gamma(1+\tau-s_1-s_4)\mathring{g}(\tau,s_1)d\tau.$$

To proceed, we focus on the context of the minus sign. In this case, the bound shown in Lemma 3.2 plays a crucial rôle. By definition, the integral transform in question is written as

$$\mathscr{L}^-\Psi_{\mathsf{s}}^-(t) = 8i^{-\kappa}\cosh(\pi t)\int_0^\infty\Psi_{\mathsf{s}}^-(\eta)K_{2it}(\eta)\frac{d\eta}{\eta}.$$

Inserting equation (3.15), we obtain an absolutely convergent integral, for we have the growth condition

$$K_{2it}(\eta) \ll \begin{cases} |\log\eta| & \text{as } \eta\to 0, \\ \exp(-\eta) & \text{as } \eta\to\infty. \end{cases}$$

Substituting Heaviside's integral formula

$$\int_0^\infty K_{2v}(\eta)\left(\frac{\eta}{2}\right)^{2s-1}d\eta = \frac{\Gamma(s+v)\Gamma(s-v)}{2} \quad \text{for} \quad \Re(s) > |\Re(v)|, \tag{3.25}$$

we are left with the expression

$$\mathscr{L}^-\Psi_{\mathsf{s}}^-(t) = -i^{-\kappa}\cosh(\pi t)\int_{(B)}\cos\left(\pi\left(s_1+\frac{s_3+s_4}{2}-\tau\right)\right)\Gamma\left(\frac{s_1+s_2+s_3+s_4-1}{2}+it-\tau\right)$$
$$\times\Gamma\left(\frac{s_1+s_2+s_3+s_4-1}{2}-it-\tau\right)\Gamma(1+\tau-s_1-s_3)\Gamma(1+\tau-s_1-s_4)\mathring{g}(\tau,s_1)\frac{d\tau}{\pi i}.$$

In anticipation of future simplifications of $\mathscr{L}^\pm\Psi_{\mathsf{s}}^\pm(t)$, we define the the integral transforms

$$\Phi_{\mathsf{s}}^+(\xi) := -2i(2\pi)^{s_1-s_2-2}(i\xi\cot(\pi\xi))^{\kappa}\cos\left(\frac{\pi(s_3-s_4)}{2}\right)\int_{-i\infty}^{i\infty}\sin\left(\frac{\pi(s_1+s_2+s_3+s_4-\kappa-2\tau)}{2}\right)$$
$$\times\Gamma\left(\frac{s_1+s_2+s_3+s_4-1}{2}+\xi-\tau\right)\Gamma\left(\frac{s_1+s_2+s_3+s_4-1}{2}-\xi-\tau\right)$$
$$\times\Gamma(1+\tau-s_1-s_3)\Gamma(1+\tau-s_1-s_4)\mathring{g}(\tau,s_1)d\tau$$

$$\tag{3.26}$$

---

[5]Motohashi [55, Page 161] made a minor mistake in the argument of a gamma function.

and

$$\Phi_{\mathsf{s}}^{-}(\xi) := 2i^{1-\kappa}(2\pi)^{s_1-s_2-2}\cos(\pi\xi)\int_{-i\infty}^{i\infty}\cos\left(\pi\left(s_1+\frac{s_3+s_4}{2}-\tau\right)\right)\Gamma\left(\frac{s_1+s_2+s_3+s_4-1}{2}+\xi-\tau\right)$$
$$\times\,\Gamma\left(\frac{s_1+s_2+s_3+s_4-1}{2}-\xi-\tau\right)\Gamma(1+\tau-s_1-s_3)\Gamma(1+\tau-s_1-s_4)\mathring{g}(\tau,s_1)d\tau$$

(3.27)

along with the original integral transform

$$\Xi_{\mathsf{s}}(\xi) := \int_{-i\infty}^{i\infty}\frac{\Gamma\left(\xi+\frac{s_1+s_2+s_3+s_4-1}{2}-\tau\right)}{\Gamma\left(\xi+\frac{3-s_1-s_2-s_3-s_4}{2}+\tau\right)}\Gamma(1+\tau-s_1-s_3)\Gamma(1+\tau-s_1-s_4)\mathring{g}(\tau,s_1)\frac{d\tau}{2\pi i}. \quad (3.28)$$

Here the contour in the formula for $\Phi_{\mathsf{s}}^{+}(\xi)$ is curved to ensure that the poles of the first two gamma factors in the integrand lie to the right of the contour and those of the other factors are on the left of the contour. In addition, the variables $\xi, s_1, s_2, s_3, s_4$ are assumed to be such that the path can be drawn. The contour in the definition of $\Phi_{\mathsf{s}}^{-}(\xi)$ is chosen in just the same way. On the other hand, the contour in $\Xi$ separates the poles of $\Gamma(\xi+(s_1+s_2+s_3+s_4-1)/2-\tau)$ and those of $\Gamma(1+\tau-s_1-s_3)\Gamma(1+\tau-s_1-s_4)\mathring{g}(\tau,s_1)$ to the left and the right of the contour. An explicit formulation of the integral transforms can be executed:

**Proposition 3.5.** *With the notation above, we have that*

$$\Phi_{\mathsf{s}}^{+}(\xi) = -(-i\xi)^{\kappa}\frac{(2\pi)^{s_1-s_2}}{2\sin(\pi\xi)}\cos\left(\frac{\pi(s_3-s_4)}{2}\right)\left\{\Xi_{\mathsf{s}}(\xi)-(-1)^{\kappa}\Xi_{\mathsf{s}}(-\xi)\right\},$$
$$\Phi_{\mathsf{s}}^{-}(\xi) = i^{-\kappa}\frac{(2\pi)^{s_1-s_2}}{2\sin(\pi\xi)}\left\{\sin\left(\pi\left(\frac{s_2-s_1}{2}+\xi\right)\right)\Xi_{\mathsf{s}}(\xi)-\sin\left(\pi\left(\frac{s_2-s_1}{2}-\xi\right)\right)\Xi_{\mathsf{s}}(-\xi)\right\},$$

*provided the left-hand sides are well-defined. We also have for real $t$ and $\mathsf{s} = (s_1, s_2, s_3, s_4) \in \mathcal{E}_4$ that*

$$\mathscr{L}^{+}\Psi_{\mathsf{s}}^{+}(t) = (2\pi)^{1-s_1+s_2}\Phi_{\mathsf{s}}^{+}(it),$$
$$\mathscr{L}^{-}\Psi_{\mathsf{s}}^{-}(t) = (2\pi)^{1-s_1+s_2}\Phi_{\mathsf{s}}^{-}(it).$$

*For integral $k \equiv \kappa \pmod 2$ and $(s_1, s_2, s_3, s_4) \in \mathcal{E}_4$, we have that*

$$\mathscr{L}^{\mathrm{hol}}\Psi_{\mathsf{s}}^{+}(k) = 2\pi i^{k}\cos\left(\frac{\pi(s_3-s_4)}{2}\right)\Xi_{\mathsf{s}}\left(\frac{k-1}{2}\right).$$

*Proof.* The proof proceeds as in [55, Lemma 4.4] *mutatis mutandis.* □

**Remark 3.6.** From Proposition 3.3, we observe that all Dirichlet characters in question are of the form $\psi^2$. It therefore suffices to treat the case of $\kappa(\psi^2) = 0$ in what follows. This should not cause any problem, although we have suppressed the dependence on $\psi$ from the functions $\Phi_{\mathsf{s}}^{\pm}(\xi)$ and $\Xi_{\mathsf{s}}(\xi)$.

### 3.5. Deduction of the cubic moment side

We devote this subsection to the derivation of the cubic moment side. The absolute convergence that we would like to check is apparent as far as the double summation over the variables $m, n$ is concerned, since we have the bound in equation (2.22) on the Hecke eigenvalues and $\mathsf{s} = (s_1, s_2, s_3, s_4) \in \mathcal{E}_4$. Hence, the chief issue is reduced to bounding $\mathscr{L}^{\pm}\Psi_{\mathsf{s}}^{\pm}$, and then Proposition 3.5 reduces our task to the analysis of the function $\Xi$. If real $t$ and positive integral $k$ tend to infinity, we have uniformly for any compact subset of $\mathcal{E}_4$ that

$$\Xi_s(it) \ll |t|^{-A}, \qquad \Xi_s\left(\frac{k-1}{2}\right) \ll k^{-A}.$$

### 3.5.1. Computation of $\mathcal{J}_\pm$

We are able to insert equation (3.21) into equation (3.19) and change the order of summations and integrals as long as we work inside $\mathcal{E}_4$. The replacement of $\psi$ (in the discussion of Section 3.4) with $\psi^2$ is necessary. We now evaluate Rankin–Selberg $L$-functions involving the divisor function defined as

$$\mathcal{T}(s,u;f;\psi) := \sum_{m=1}^\infty \frac{\overline{\psi}(m)\sigma_u(m)\lambda_f(m)}{m^s}, \qquad \mathcal{T}(s,u;E;\psi) := \sum_{m=1}^\infty \frac{\overline{\psi}(m)\sigma_u(m)\lambda_f(m,t)}{m^s}.$$

One then establishes the following lemma:

**Lemma 3.7.** *For $f \in \mathcal{B}_\kappa(dq,\psi^2)$ and an Eisenstein series $E(z,1/2+it,f)$ with $f \in \mathcal{B}(\psi_1,\psi_2)$, we have*

$$\mathcal{T}(s,u;f;\psi) = \frac{L(s, f \otimes \overline{\psi})L(s-u, f \otimes \overline{\psi})}{\zeta^q(2s-u)},$$

$$\mathcal{T}(s,u;E;\psi) = \frac{L(s+it, \overline{\psi}\psi_2)L(s-it, \overline{\psi}\psi_1)L(s-u+it, \overline{\psi}\psi_2)L(s-u-it, \overline{\psi}\psi_1)}{\zeta^q(2s-u)}.$$

**Remark 3.8.** The level $dq$ in Lemma 3.7 can be replaced with any positive integer that is divided by the modulus $q$ of the central character.

*Proof.* For the first assertion, we use the multiplicativity relation in equation (2.21) for the Hecke eigenvalues, getting

$$\mathcal{T}(s,u;f;\psi) = \sum_{m=1}^\infty \sum_{n=1}^\infty \frac{\overline{\psi}(mn)\lambda_f(mn)}{m^s n^{s-u}} = \sum_{m=1}^\infty \sum_{n=1}^\infty \frac{\overline{\psi}(mn)}{m^s n^{s-u}} \sum_{d|(m,n)} \mu(d)\psi(d)^2 \lambda_f\left(\frac{m}{d}\right)\lambda_f\left(\frac{n}{d}\right).$$

Upon exchanging the order of summations, the right-hand side equals

$$\sum_{(d,q)=1} \frac{\mu(d)}{d^{2s-u}} \sum_{m=1}^\infty \frac{\overline{\psi}(m)\lambda_f(m)}{m^s} \sum_{n=1}^\infty \frac{\overline{\psi}(n)\lambda_f(n)}{n^{s-u}} = \frac{L(s, f \otimes \overline{\psi})L(s-w, f \otimes \overline{\psi})}{\zeta^q(2s-u)}.$$

In order to establish the second assertion, we make use of Lemma 2.5 to derive

$$\begin{aligned}
\mathcal{T}(s,u;E;\psi) &= \sum_{m=1}^\infty \frac{\overline{\psi}(m)\sigma_u(m)}{m^s} \sum_{a|m} \psi_1(a)\psi_2\left(\frac{m}{a}\right)\left(\frac{a^2}{m}\right)^{it} \\
&= \sum_{(c,q)=1} \frac{\mu(c)}{c^{2s-u}} \sum_{a=1}^\infty \frac{\overline{\psi}(a)\psi_1(a)\sigma_w(a)}{a^{s-it}} \sum_{m=1}^\infty \frac{\overline{\psi}(m)\psi_2(m)\sigma_u(m)}{m^{s+it}} \\
&= \frac{L(s+it, \overline{\psi}\psi_2)L(s-it, \overline{\psi}\psi_1)L(s-u+it, \overline{\psi}\psi_2)L(s-u-it, \overline{\psi}\psi_1)}{\zeta^q(2s-u)}.
\end{aligned}$$

This concludes the proof of Lemma 3.7. □

We proceed to compute $\mathcal{A}_{\infty\infty 0}^{\mathrm{Maaß}}(m, \pm n; \mathscr{L}^{\pm}\Psi_{\mathsf{s}}^{\pm}; \psi^2)$ by virtue of Lemma 3.7, obtaining

$$\sum_{n=1}^{\infty} \frac{\overline{\psi}(n)\sigma_{s_4-s_3}(n)\rho_f(\pm n)}{n^{(1-s_1+s_2-s_3+s_4)/2}} = \epsilon_f^{(1\mp 1)/2}\rho_f(1)\frac{L\left(\frac{1-s_1+s_2+s_3-s_4}{2}, f\otimes\overline{\psi}\right)L\left(\frac{1-s_1+s_2-s_3+s_4}{2}, f\otimes\overline{\psi}\right)}{\zeta^q(1-s_1+s_2)},$$

where $f \in \mathcal{B}_{\kappa}(dq, \psi^2)$. There is beauty in the above quotient since we have $\zeta^q(1-s_1+s_2)^{-1}$, which completely cancels out with the term $\zeta^q(1-s_1+s_2)$ appearing in the expression for $\mathcal{J}_{\pm}$ in Proposition 3.3. One can calculate the Dirichlet series for the sum over $m$ in a similar manner, obtaining

$$\sum_{m=1}^{\infty} \frac{\psi(m)\overline{\rho_f(m)}}{m^{(s_1+s_2+s_3+s_4-1)/2}} = \overline{\rho_f(1)}L\left(\frac{s_1+s_2+s_3+s_4-1}{2}, f\otimes\overline{\psi}\right),$$

where $f \in \mathcal{B}_{\kappa}(dq, \psi^2)$. We obtain the central $L$-value $L(1/2, f\otimes\overline{\psi})^3$ if we take the limit $(s_1, s_2, s_3, s_4) \mapsto (1/2, 1/2, 1/2, 1/2)$. A similar procedure can be applied to the Eisenstein contribution to deduce the six Dirichlet $L$-functions in equation (1.4). The eventual form of the cubic moment side looks like

$$\mathcal{J}_{\pm} = \mathcal{J}_{\pm}^{\mathrm{Maaß}} + \mathcal{J}_{\pm}^{\mathrm{Eis}} + \delta_{\pm=+}\mathcal{J}_{+}^{\mathrm{hol}}$$

with

$$\mathcal{J}_{\pm}^{\mathrm{Maaß}} := \frac{q^{s_2-s_1-2}}{\tau(\overline{\chi})} \sum_{d|q} \frac{\mu(d)}{d} \sum_{\psi \,(\mathrm{mod}\, q)}^{\#} \mathcal{H}_{\pm}(\chi, \psi) \sum_{f \in \mathcal{B}_{\kappa}^*(dq, \psi^2)} \epsilon_f^{(1\mp 1)/2}$$

$$\times \frac{L\left(\frac{1-s_1+s_2+s_3-s_4}{2}, f\otimes\overline{\psi}\right)L\left(\frac{1+s_2-s_3+s_4}{2}, f\otimes\overline{\psi}\right)L\left(\frac{s_1+s_2+s_3+s_4-1}{2}, f\otimes\overline{\psi}\right)}{L(1, \mathrm{Ad}^2 f)}\Phi_{\mathsf{s}}^{\pm}(it_f),$$

$$\mathcal{J}_{\pm}^{\mathrm{Eis}} := \frac{2q^{s_2-s_1-1}}{\tau(\overline{\chi})\varphi(q)} \sum_{d|q} \frac{\mu(d)}{d} \sum_{\psi \,(\mathrm{mod}\, q)}^{\#} \mathcal{H}_{\pm}(\chi, \psi) \sum_{\psi_1\psi_2=\psi^2} \sum_{f \in \mathcal{B}(\psi_1, \psi_2)} \int_{-\infty}^{\infty} \frac{\mathcal{S}_f(t; s_1, s_2, s_3, s_4)}{|L(1+2it, \overline{\psi_1}\psi_2)|^2}\Phi_{\mathsf{s}}^{\pm}(it)\frac{dt}{2\pi},$$

$$\mathcal{J}_{+}^{\mathrm{hol}} := \frac{q^{-2}}{\tau(\overline{\chi})}\left(\frac{2\pi}{q}\right)^{s_1-s_2}\cos\left(\frac{\pi(s_3-s_4)}{2}\right) \sum_{d|q} \frac{\mu(d)}{d} \sum_{\psi \,(\mathrm{mod}\, q)}^{\#} \mathcal{H}_{+}(\chi, \psi) \sum_{\substack{k>\kappa \\ k\equiv\kappa\,(\mathrm{mod}\,2)}} \sum_{f \in \mathcal{B}_k^*(dq, \psi^2)} i^k$$

$$\times \frac{L\left(\frac{1-s_1+s_2+s_3-s_4}{2}, f\otimes\overline{\psi}\right)L\left(\frac{1-s_1+s_2-s_3+s_4}{2}, f\otimes\overline{\psi}\right)L\left(\frac{s_1+s_2+s_3+s_4-1}{2}, f\otimes\overline{\psi}\right)}{L(1, \mathrm{Ad}^2 f)}\Xi_{\mathsf{s}}\left(\frac{k-1}{2}\right).$$

Here we set $\mathcal{H}_{\pm}(\chi, \psi) = \psi(\mp 1)\tau(\psi)\tau(\overline{\chi\psi})J(\psi, \psi)$ and used the identity $\varphi(dq) = d\varphi(q)$ for $q$ prime.

### 3.5.2. Computation of $\mathcal{E}_{\pm}$

We initially handle $\mathcal{A}_{\infty\infty 0}^{\mathrm{Maaß}}(m, \pm n_0 n; \mathscr{L}^{\pm}\Psi_{\mathsf{s}}^{\pm}; \psi_0)$. One can rewrite the direct sum decomposition in equation (2.20) in order to imitate the formulation of Blomer–Milićević [14], namely

$$\mathcal{B}_0(\Gamma_0(r^2)) = \bigsqcup_{r_1 r_2 = r^2} \bigsqcup_{f \in \mathcal{B}_0^*(\Gamma_0(r_1))} \bigsqcup_{d|r_2} \iota_d f \cdot \mathbb{C},$$

where the first two sums are orthogonal, but the last one is not orthogonal and needs to be orthogonalised by Gram–Schmidt. Here we also define $(\iota_d f)(z) := f(dz)$. By [14, Lemma 9], the set of functions

$$\left\{f^{(g)} := \sum_{d|g} \xi_g(d)(\iota_d f) : g \mid r_2\right\}$$

is an orthonormal basis of the space $\bigsqcup_{d|r_2} \iota_d f \cdot \mathbb{C}$ with $f \in \mathcal{B}_0^*(\Gamma_0(r_1))$. Then the Fourier coefficients are

$$\rho_{f^{(g)}\mathfrak{a}}(n) = \sum_{d|g} \xi_g(d) d^{1/2} \rho_{f\mathfrak{a}}\left(\frac{n}{d}\right)$$

with the convention that $\rho_{f\mathfrak{a}}(x) = 0$ for $x \notin \mathbb{Z}$. We calculate the sum over $n$ in equation (3.20) as

$$\sum_{n=1}^{\infty} \frac{\psi(n)\sigma_{s_4-s_3}(n)\rho_{f0}(\pm n/d)}{n^{(1-s_1-s_2+s_4)/2}} = \delta_{d=1}\epsilon_f^{(1\mp1)/2}\omega_f\rho_f(1)\frac{L\left(\frac{1-s_1+s_2+s_3-s_4}{2}, f\otimes\psi\right)L\left(\frac{1-s_1+s_2-s_3+s_4}{2}, f\otimes\psi\right)}{\zeta^r(1-s_1+s_2)},$$

where $\omega_f$ stands for the root number for a newform $f \in \mathcal{B}_0^*(\Gamma_0(r_1))$. We denote by $\psi$ a Dirichlet character modulo $r$ induced by the primitive character $\psi^*$ modulo $r^*$. We regard the character modulo 1 as primitive. Since the central character is trivial and the Hecke eigenvalues are real (Remark 2.10), we use Lemma 2.1 to infer that

$$\sum_{m=1}^{\infty} \frac{\tau(\psi, m)\overline{\rho_f(m/d)}}{m^{(s_1+s_2+s_3+s_4-1)/2}} = \frac{\overline{\rho_f(1)}\tau(\psi^*)\psi^*(d)}{d^{(s_1+s_2+s_3+s_4-1)/2}}\Sigma_{r/r^*}(s_1, s_2, s_3, s_4; \psi)L\left(\frac{s_1+s_2+s_3+s_4-1}{2}, f\otimes\psi^*\right),$$

where $\Sigma_q(s_1, s_2, s_3, s_4; \psi)$ is the multiplicative function

$$\sum_{c|q} \mu\left(\frac{q}{c}\right)\psi^*\left(\frac{q}{c}\right)c^{(3-s_1-s_2-s_3-s_4)/2}\sum_{e|c}\mu(e)\psi_0(e)\psi^*(e)e^{(1-s_1-s_2-s_3-s_4)/2}\lambda_f\left(\frac{c}{e}\right)$$

and $\psi_0$ is the principal character modulo $r_1$. On the other hand, we obtain via elementary calculations that

$$\sum_{n|r^\infty} \frac{\tau(\psi, n)\sigma_w(n)\lambda_f(n)}{n^s} = \tau(\psi^*)\sum_{c|r/r^*}c\psi^*\left(\frac{r}{cr^*}\right)\mu\left(\frac{r}{cr^*}\right)\sum_{n|r^\infty}\frac{\psi^*(n)\sigma_w(nc)\lambda_f(nc)}{(nc)^s}, \qquad (3.29)$$

where $s = (1-s_1+s_2-s_3+s_4)/2$ and $w = s_4 - s_3$. We hence observe that the left-hand side of equation (3.29) equals $\tau(\psi^*)$ when $\psi$ is the quadratic character modulo $q$ or the character modulo 1. The problem occurs if $\psi$ is the principal character modulo a prime $q$. In this case, a brute force computation gives

$$\sum_{c|q} c^{1-s}\mu\left(\frac{q}{c}\right)\sum_{n|q^\infty}\frac{\sigma_w(nc)\lambda_f(nc)}{n^s} = \varphi(q)\left(1-\frac{\lambda_f(q)}{q^s}\right)^{-1}\left(1-\frac{\lambda_f(q)}{q^{s-w}}\right)^{-1} - q.$$

We then conclude that

$$\sum_{n|r^\infty} \frac{\tau(\psi, n)\sigma_w(n)\lambda_f(n)}{n^{(1-s_1+s_2-s_3+s_4)/2}} = \tau(\psi^*)\Pi_{r/r^*}(s_1, s_2, s_3, s_4; \psi),$$

where $\Pi_q(s_1, s_2, s_3, s_4; \psi)$ is the multiplicative function

$$q\sum_{c|q}\mu\left(\frac{q}{c}\right)\prod_{p|c}\left(1-\frac{\psi^*(p)}{p}\right)\left(1-\frac{\psi^*(p)\lambda_f(p)}{p^{(1-s_1+s_2+s_3-s_4)/2}}\right)^{-1}\left(1-\frac{\psi^*(p)\lambda_f(p)}{p^{(1-s_1+s_2-s_3+s_4)/2}}\right)^{-1}.$$

This reasoning proceeds verbatim for the Eisenstein and holomorphic spectra, giving similar expressions (see [70] for an extension of Weisinger's theory and newform Eisenstein series). Thus we are left with

$$\mathcal{E}_\pm = \mathcal{E}_\pm^{\text{Maaß}} + \mathcal{E}_\pm^{\text{Eis}} + \delta_{\pm=+}\mathcal{E}_+^{\text{hol}},$$

where

$$\mathcal{E}_\pm^{\text{Maaß}} := \frac{1}{\tau(\overline{\chi})q} \sum_{r|q} \frac{r^{1-s_1+s_2}}{\varphi(r)} \sum_{\psi(\text{mod } r)}^\flat \mathcal{G}_\pm(\chi,\psi) \sum_{r_1 r_2 = r^2} \sum_{f \in \mathcal{B}_0^*(\Gamma_0(r_1))} \epsilon_f^{(1\mp1)/2} \omega_f \, \Omega_{r_1,r_2}(s_1,s_2,s_3,s_4;f;\psi)$$

$$\times \frac{L\left(\frac{1-s_1+s_2+s_3-s_4}{2}, f\otimes\psi\right) L\left(\frac{1-s_1+s_2-s_3+s_4}{2}, f\otimes\psi\right) L\left(\frac{s_1+s_2+s_3+s_4-1}{2}, f\otimes\psi^*\right)}{L(1,\text{Ad}^2 f)} \Phi_{\mathsf{s}}^\pm(it_f),$$

$$\mathcal{E}_\pm^{\text{Eis}} := \frac{2}{\tau(\overline{\chi})q} \sum_{r|q} \frac{r^{1-s_1+s_2}}{\varphi(r)} \sum_{\psi(\text{mod } r)}^\flat \mathcal{G}_\pm(\chi,\psi) \sum_{r_1 r_2 = r^2} \frac{r_1^2}{\varphi(r_1)} \sum_{\psi_1\psi_2=\psi_0}^* \sum_{f\in\mathcal{B}(\psi_1,\psi_2)}$$

$$\times \Omega_{r_1,r_2}(s_1,s_2,s_3,s_4;f;\psi) \int_{-\infty}^\infty \overline{\rho_f(1,t)}\rho_{f0}(1,t)\mathcal{S}_f(t;s_1,s_2,s_3,s_4)\Phi_{\mathsf{s}}^\pm(it)\frac{dt}{2\pi},$$

$$\mathcal{E}_+^{\text{hol}} := \frac{(2\pi)^{s_1-s_2}}{\tau(\overline{\chi})q} \cos\left(\frac{\pi(s_3-s_4)}{2}\right) \sum_{r|q} \frac{r^{1-s_1+s_2}}{\varphi(r)} \sum_{\psi(\text{mod } r)}^\flat \mathcal{G}_+(\chi,\psi)$$

$$\times \sum_{r_1 r_2 = r^2} \sum_{\substack{k > \kappa \\ k \equiv \kappa(\text{mod } 2)}} \sum_{f\in\mathcal{B}_k^*(\Gamma_0(r_1))} i^k \omega_f \, \Omega_{r_1,r_2}(s_1,s_2,s_3,s_4;f;\psi)$$

$$\times \frac{L\left(\frac{1-s_1+s_2+s_3-s_4}{2}, f\otimes\psi\right) L\left(\frac{1-s_1+s_2-s_3+s_4}{2}, f\otimes\psi\right) L\left(\frac{s_1+s_2+s_3+s_4-1}{2}, f\otimes\psi^*\right)}{L(1,\text{Ad}^2 f)} \Xi_{\mathsf{s}}\left(\frac{k-1}{2}\right),$$

where we define $\mathcal{G}_\pm(\chi,\psi) = \psi(\mp1)\tau(\psi)\tau(\psi^*)^2\tau(\overline{\chi}\psi)$ and

$$\Omega_{r_1,r_2}(s_1,s_2,s_3,s_4;f;\psi) = \frac{\varphi(r_1)}{r_1^2} \sum_{g|r_2} \xi_g(1) \sum_{d|g} \xi_g(d)\psi^*(d)d^{(2-s_1-s_2-s_3-s_4)/2}$$

$$\times \Sigma_{r/r^*}(s_1,s_2,s_3,s_4;\psi)\Pi_{r/r^*}(s_1,s_2,s_3,s_4;\psi).$$

These terms are complicated, but they are the same size as the main contributions $\mathcal{J}_\pm$ in estimations. Thus one can ignore them in the process of using Motohashi's formula to deduce some moment bounds.

### 3.6. Endgame: analytic continuation

The aim of this subsection is to show that the spectral decomposition obtained in the preceding subsection can be continued to the whole complex plane $\mathbb{C}^4$ and hence to conclude the proof of Theorem 1.1. The main step is identical to that of Motohashi's monograph [55], and we stress the following lemmata:

**Lemma 3.9** (Motohashi [55, Lemma 4.7]). *The function $\Xi_{\mathsf{s}}(\xi)$ is meromorphic in the domain*

$$\mathcal{B}_4^* = \{\xi : \Re(\xi) > -cA\} \times \mathcal{B}_4$$

*for a fixed small constant $c > 0$ and holomorphic in $\mathcal{B}_4^* \setminus \mathcal{N}$, where*

$$\mathcal{N} = \left\{ (\xi,s_1,s_2,s_3,s_4) \ : \ \begin{array}{l} \text{at least one of } \xi + \frac{s_1+s_2+s_3+s_4-1}{2}, \xi + \frac{1-s_1+s_2+s_3-s_4}{2}, \\ \xi + \frac{1-s_1+s_2-s_3+s_4}{2} \text{ equals a non-positive integer} \end{array} \right\}.$$

*Moreover, if $|\xi|$ tends to infinity in any fixed vertical or horizontal strips while satisfying $\Re(\xi) > -cA$, then we have uniformly in $\mathcal{B}_4$ that*

$$\Xi_{\mathsf{s}}(\xi) \ll |\xi|^{-cA}. \tag{3.30}$$

**Lemma 3.10** (Motohashi [55], Lemma 4.8]). *If $(\xi, s_1, s_2, s_3, s_4)$ is such that the path in equation (3.28) can be drawn in a vertical strip contained in the half plane $\Re(\tau) > 0$, then we have that*

$$\Xi_s(\xi) = \frac{\Gamma(\alpha)\Gamma(\beta)}{\Gamma(\gamma)} \int_0^\infty y^{\xi + (s_1 + s_2 + s_3 + s_4 - 3)/2} G(y, s_1) \, {}_2F_1(\alpha, \beta; \gamma; -y) \, dy,$$

*where $G(y, s)$ is defined in equation (3.4) and ${}_2F_1(\alpha, \beta; \gamma; y)$ is the hypergeometric function with*

$$\alpha = \xi + \frac{1 - s_1 + s_2 + s_3 - s_4}{2}, \qquad \beta = \xi + \frac{1 - s_1 + s_2 - s_3 + s_4}{2}, \qquad \gamma = 1 + 2\xi.$$

We assume that $(s_1, s_2, s_3, s_4) \in \mathcal{B}_4$ and examine the consequence of Lemma 3.9 for the contribution $\mathcal{J}_+^{\text{Maaß}}$ from Maaß forms. If $t_f \geqslant 3B$, then $\Xi_s(\pm i t_f)$ is holomorphic and $O(t_f^{-cA})$ uniformly in $\mathcal{B}_4$ with an absolute constant $c$. Hence via Proposition 3.5, the function $\mathscr{L}^+ \Psi_s^+(t_f)$ is holomorphic and of exponential decay with respect to $t_f$ uniformly in $\mathcal{B}_4$. This indicates that $\mathcal{J}_+^{\text{Maaß}}$ exists as a meromorphic function inside $\mathcal{B}_4$. The same observation holds for $\mathcal{J}_+^{\text{hol}}$. As for the function $\mathcal{J}_-^{\text{Maaß}}$, we essentially need equation (3.30). From Proposition 3.5, we have $\mathscr{L}^- \Psi_s^-(t_f) = O(t_f^{-cA})$ uniformly in $\mathcal{B}_4$ provided $t_f \geqslant 3B$. Hence $\mathcal{J}_-^{\text{Maaß}}$ is meromorphic inside $\mathcal{B}_4$. This discussion works for $\mathcal{E}_\pm^{\text{Maaß}}$ and $\mathcal{E}_+^{\text{hol}}$ as well.

Thus it remains to contemplate $\mathcal{J}_\pm^{\text{Eis}}$ and $\mathcal{E}_\pm^{\text{Eis}}$. To this end, we assume first that $s = (s_1, s_2, s_3, s_4) \in \mathcal{E}_4$ and set

$$\mathcal{J}^{\text{Eis}}(s; g; \chi) = \mathcal{J}_+^{\text{Eis}}(s; g; \chi) + \mathcal{J}_-^{\text{Eis}}(s; g; \chi).$$

Using Proposition 3.5, one derives

$$\Phi_s^+(it) = -\frac{(2\pi)^{s_1 - s_2}}{2i \sinh(\pi t)} \cos\left(\frac{\pi(s_3 - s_4)}{2}\right) \left\{ \Xi_s(it) - \Xi_s(-it) \right\},$$

$$\Phi_s^-(it) = \frac{(2\pi)^{s_1 - s_2}}{2i \sinh(\pi t)} \left\{ \sin\left(\pi\left(\frac{s_2 - s_1}{2} + it\right)\right) \Xi_s(it) - \sin\left(\pi\left(\frac{s_2 - s_1}{2} - it\right)\right) \Xi_s(-it) \right\}.$$

This yields an expression of $\mathcal{J}^{\text{Eis}}(s; g; \chi)$ in terms of $\Xi_s(it)$. We need to shift the contour in the same way as on pages 170–171 of Motohashi's monograph [55]. This is because $\mathcal{S}_f(t; s_1, s_2, s_3, s_4)$ has poles when $\psi_1 = \psi_2 = \psi$. The same procedure applies to the term $\mathcal{E}_\pm^{\text{Eis}}$ coming from the contribution of the principal character and quadratic character. Gathering these observations, we have proven the following lemma:

**Lemma 3.11.** *The function $\mathcal{OD}_4^\dagger$ continues meromorphically to $\mathcal{B}_4$, and the decomposition in equation (3.2) holds throughout this domain. In particular, the formula in equation (1.5) is valid for $(s_1, s_2, s_3, s_4) = (1/2, 1/2, 1/2, 1/2)$.*

This concludes the proof of Theorem 1.1.

## 4. Implications of Theorem 1.1

In this section, we aim to establish a $q$-aspect variant of Iwaniec's short interval fourth moment bound and the twelfth moment bound for Dirichlet $L$-functions without an average over Dirichlet characters.

### 4.1. Proofs of Corollaries 1.3 and 1.4

It suffices to specialise the test function in Theorem 1.1 as

$$g(t) = \frac{1}{\sqrt{\pi} H} \exp\left(-\left(\frac{t - T}{H}\right)^2\right) \tag{4.1}$$

with the parameter $H$ at one's disposal. We assume that

$$T^{1/2} \leqslant H \leqslant T(\log T)^{-1}.$$

Motohashi [55, (5.1.40)–(5.1.42)] evaluated the corresponding integral transform $\Xi$, and his expression can be used in our context. We then derive an asymptotic formula almost identical to [55, Theorem 5.1]. Upon estimating the spectral sum by absolute values, an oscillatory component in the summand evanishes, and we obtain the bound

$$\int_T^{T+H} \left| L\left(\frac{1}{2} + it, \chi\right) \right|^4 dt \ll_\epsilon H^{1+\epsilon} q^\epsilon + \frac{H^{3/2}}{T} q^{-2+\epsilon} \sum_{\psi \,(\mathrm{mod}\ q)} \sum_{\substack{f \in \mathcal{B}_\kappa^*(q^2, \psi^2) \\ t_f \ll T/H}} L\left(\frac{1}{2}, f \otimes \overline{\psi}\right)^3. \qquad (4.2)$$

The truncation of the sum is justified since the Taylor expansion of $\exp(-(Ht_f/T)^2/4)$ in [55, (5.1.44)] implies that the rest of the sum is smaller than the main term in equation (4.2). Notice that the central $L$-values in equation (4.2) are nonnegative, and it thus follows from the work of Petrow–Young [63, Theorems 1.2 & 1.3] that

$$H^{1+\epsilon} q^\epsilon + \frac{H^{3/2}}{T} q^{-1+\epsilon} \left(\frac{qT}{H}\right)^{2+\epsilon} \ll H^{1+\epsilon} q^\epsilon + \left(\frac{qT}{\sqrt{H}}\right)^{1+\epsilon}.$$

Optimising $H = (qT)^{2/3}$ concludes the proof of Corollaries 1.3 and 1.4.

### 4.2. Proofs of Corollaries 1.5 and 1.6

This section is devoted to proving the twelfth moment bound. We follow the argument of [34] verbatim, obtaining the expression[6]

$$\sum_{r=1}^R \int_{t_r}^{t_r+T_0} \left| L\left(\frac{1}{2} + it, \chi\right) \right|^4 dt \ll_\epsilon T_0 q^{-2+\epsilon} \sum_{\psi \,(\mathrm{mod}\ q)} \sum_{\substack{f \in \mathcal{B}_\kappa^*(q^2, \psi^2) \\ t_f \leqslant T T_0^{-1} \sqrt{\log T}}} t_f^{-1/2} \frac{L(1/2, f \otimes \overline{\psi})^3}{L(1, \mathrm{Ad}^2 f)}$$

$$\times \sum_{r=1}^R t_r^{-1/2} \exp\left(-\left(\frac{T_0 t_f}{2t_r}\right)^2\right) \sin\left(t_f \log \frac{t_f}{4et_r}\right) + RT_0(qT)^\epsilon.$$

For technical reasons, it is convenient to remove the exponential factor in the last sum over $r$ by partial summation. In doing this, we obtain an error term of $O_\epsilon(RT^\epsilon T_0^{-1})$. Then we must majorise

$$\sum := \sum_{\substack{K = 2^{-m} T T_0^{-1} \sqrt{\log T} \\ m = 1, 2, \ldots}} S_K,$$

where

$$S_K := \sum_{\psi \,(\mathrm{mod}\ q)} \sum_{\substack{f \in \mathcal{B}_\kappa^*(q^2, \psi^2) \\ K \leqslant t_f \leqslant 2K}} t_f^{-1/2} \frac{L(1/2, f \otimes \overline{\psi})}{L(1, \mathrm{Ad}^2 f)} \left| \sum_{r=1}^R t_r^{-1/2 - it_f} \right|.$$

---

[6]The variables $t_r$ and spectral parameters $t_f$ should not be confused.

The Cauchy–Schwarz inequality leads one to

$$S_K \leqslant \left( \sum_{\psi \,(\mathrm{mod}\, q)} \sum_{\substack{f \in \mathcal{B}_\kappa^*(q^2, \psi^2) \\ K \leqslant t_f \leqslant 2K}} t_f^{-1} \frac{L(1/2, f \otimes \overline{\psi})^4}{L(1, \mathrm{Ad}^2 f)} \right)^{1/2}$$

$$\times \left( \sum_{\psi \,(\mathrm{mod}\, q)} \sum_{\substack{f \in \mathcal{B}_\kappa^*(q^2, \psi^2) \\ K \leqslant t_f \leqslant 2K}} \frac{L(1/2, f \otimes \overline{\psi})^2}{L(1, \mathrm{Ad}^2 f)} \left| \sum_{r=1}^{R} t_r^{-1/2 - it_f} \right|^2 \right)^{1/2}.$$

Applying the bound of Petrow–Young [63, Theorem 7.6], the fourth moment is bounded as

$$\sum_{\psi \,(\mathrm{mod}\, q)} \sum_{\substack{f \in \mathcal{B}_\kappa^*(q^2, \psi^2) \\ K \leqslant t_f \leqslant 2K}} t_f^{-1} \frac{L(1/2, f \otimes \overline{\psi})^4}{L(1, \mathrm{Ad}^2 f)} \ll_\epsilon q^{3+\epsilon} K^{1+\epsilon}.$$

The spectral large sieve in the form involving the square of twisted automorphic $L$-functions yields

$$\sum_{\psi \,(\mathrm{mod}\, q)} \sum_{\substack{f \in \mathcal{B}_\kappa^*(q^2, \psi^2) \\ K \leqslant t_f \leqslant 2K}} \frac{L(1/2, f \otimes \overline{\psi})^2}{L(1, \mathrm{Ad}^2 f)} \left| \sum_{r=1}^{R} t_r^{-1/2 - it_f} \right|^2 \ll_\epsilon q^{3+\epsilon} RKT^\epsilon T_0^{-1},$$

since $K \ll T T_0^{-1} \sqrt{\log T}$. It follows that

$$S_K \ll_\epsilon q^{3+\epsilon} R^{1/2} K^{1+\epsilon} T^\epsilon T_0^{-1/2}$$

and the summation over $K$ eventually gives

$$\sum_{r=1}^{R} \int_{t_r}^{t_r + T_0} \left| L\left(\frac{1}{2} + it, \chi\right) \right|^4 dt \ll_\epsilon \left( R T_0 + q T \sqrt{\frac{R}{T_0}} \right) (qT)^\epsilon.$$

This establishes Corollary 1.5. Corollary 1.6 follows from Jutila's trick described in [42, §6].

## A.  Twists of Maaß newforms

We evaluate the size of the conductor of $f \otimes \overline{\psi}$ via automorphic representation theory.

**Theorem A.1.** *Let $\psi$ be a primitive character modulo an odd squarefree $q$ and $d \mid q$. Let $f \in \mathcal{B}_\kappa^*(dq, \psi^2)$ be a Hecke–Maaß newform. Then the twist $f \otimes \overline{\psi}$ has conductor $q^2$ if $\psi$ is nonquadratic. If $\psi$ is quadratic, then $f \otimes \overline{\psi}$ has conductor dividing $q^2$, and all three possibilities occur: if $d = 1$, then the conductor is $q^2$; and if $d = q$, then*

$$\mathrm{cond}(f \otimes \overline{\psi}) = \begin{cases} 1 & \textit{if } f \textit{ corresponds to a principal series representation,} \\ q & \textit{if } f \textit{ corresponds to a special representation,} \\ q^2 & \textit{if } f \textit{ corresponds to a twist-minimal representation.} \end{cases}$$

An automorphic representation $f$ is named twist-minimal if it has minimal conductor among all their twists $f \otimes \psi$ by Dirichlet characters, namely $\mathrm{cond}(f) \leqslant \mathrm{cond}(f \otimes \psi)$. Theorem A.1 asserts in particular that $f \otimes \overline{\psi}$ has conductor $q^2$ when $f$ corresponds to a twist-minimal representation at the primes dividing $q$. This matches work of Booker–Lee–Strömbersson [16, Lemma 1.4]. The method to show Theorem A.1 features the use of automorphic representation theory. Nonetheless, one can also use the local Langlands correspondence and tend toward the Galois side. The only way the conductor of the tensor product could be smaller than $q^2$ is if the associated representation has a subrepresentation as an inertia representation isomorphic to $\psi$, in which case it has a quotient representation also isomorphic to $\overline{\psi}$ since its determinant is $\psi^2$. Consequently, $f$ is not twist-minimal if the conductor $\mathrm{cond}(f \otimes \overline{\psi})$ does not equal $q^2$, and there are no twist-minimal principal series representations in our context.

### A.1. Classification of representations

Let $\pi_v$ be a generic irreducible admissible representation of $\mathrm{GL}_2(\mathbb{Q}_v)$ with central character $\omega_v$. We recall that such representations can be classified into principal series representations, special representations or supercuspidal representations. Standard references for the properties of these representations are [20, 25], and the articles [32, 66] discuss the conductor exponents of these representations.

#### A.1.1. Principal series representations of $\mathrm{GL}_2(\mathbb{Q}_v)$
A principal series representation $\pi_v$ is unitarily induced from a representation of the Borel subgroup of $\mathrm{GL}_2(\mathbb{Q}_v)$, and these representations are indexed by two characters

$$\omega_{p,1} = \beta_{p,1} |\cdot|_v^{s_1}, \qquad \omega_{p,2} = \beta_{p,2} |\cdot|_v^{s_2}$$

of $\mathbb{Q}_v^\times$, where $\beta_{p,1}$ and $\beta_{p,2}$ are characters of $\mathcal{O}_v^\times$ and $s_1, s_2 \in \mathbb{C}$. We write

$$\pi_v \cong \omega_{p,1} \boxplus \omega_{p,2}.$$

This representation is irreducible and unitarisable if and only if either $s_1, s_2 \in i\mathbb{R}$ or $s_1 + s_2 \in i\mathbb{R}$ with $s_1 - s_2 \in (-1, 1)$ and $\beta_{p,1} = \beta_{p,2}$. The central character of $\pi_v$ is

$$\omega_{\pi_v} = \omega_{p,1}\omega_{p,2} = \beta_{p,1}\beta_{p,2} |\cdot|_v^{s_1+s_2}.$$

The conductor exponent $c(\pi_v)$ equals $c(\omega_{p,1}) + c(\omega_{p,2})$. If $\omega_v'$ is a character of $\mathbb{Q}_v^\times$, the twist of $\pi_v$ by $\omega_v'$ is the principal series representation

$$\pi_v \otimes \omega_v' \cong \omega_{p,1}\omega_v' \boxplus \omega_{p,2}\omega_v'.$$

If $\pi_v$ has trivial central character, then $\omega_{p,2} = \omega_{p,1}^{-1}$ so that $\beta_{p,2} = \beta_{p,1}^{-1}$, $s_2 = -s_1$ and $c(\pi_v) = 2c(\omega_{p,1})$.

#### A.1.2. Special representations of $\mathrm{GL}_2(\mathbb{Q}_v)$
A special representation is a twist of the Steinberg representation: it is the irreducible subrepresentation

$$\pi_v \cong \omega_v \mathrm{St}_v$$

of codimension one of the reducible principal series representation $\omega_v |\cdot|_v^{1/2} \boxplus \omega_v |\cdot|_v^{-1/2}$, where $\omega_v = \beta_v |\cdot|_v^s$ is a central character of $\mathbb{Q}_v^\times$ with $\beta_v$ a character of $\mathcal{O}_v^\times$ and $s \in \mathbb{C}$. The central character of $\pi_v$ is

$$\omega_{\pi_v} = \omega_v^2 = \beta_v^2 |\cdot|_v^{2s}.$$

The conductor exponent is

$$c(\pi_v) = \begin{cases} 1 & \text{if } c(\omega_v) = 0, \\ 2c(\omega_v) & \text{otherwise.} \end{cases}$$

If $\omega'_v$ is a character of $\mathbb{Q}_v^\times$, the twist of $\pi_v$ by $\omega'_v$ is the special representation

$$\pi_v \otimes \omega'_v = \omega_v \omega'_v \mathrm{St}_v.$$

### A.1.3. Supercuspidal representations of $\mathrm{GL}_2(\mathbb{Q}_v)$

A supercuspidal representation is the compact induction of a finite-dimensional representation $\rho_{\pi_v}$ of a maximal open subgroup $H$ of $\mathrm{GL}_2(\mathbb{Q}_v)$ such that $H$ is compact modulo the centre $Z(\mathbb{Q}_v)$. Every maximal open subgroup of $\mathrm{GL}_2(\mathbb{Q}_v)$ that is compact modulo the centre is conjugate to either $Z(\mathbb{Q}_v)\,\mathrm{GL}_2(\mathcal{O}_v)$ or $NK_0(\mathfrak{p}_v)$, the normaliser of $K_0(\mathfrak{p}_v)$ in $\mathrm{GL}_2(\mathbb{Q}_v)$, where $\mathfrak{p}_v$ is the maximal ideal of the ring of integers $\mathcal{O}_v$ of $\mathbb{Q}_v$. A supercuspidal representation $\pi_v$ is said to be of type I if $H$ is conjugate to $Z(\mathbb{Q}_v)\,\mathrm{GL}_2(\mathcal{O}_v)$ and of type II if $H$ is conjugate to $NK_0(\mathfrak{p}_v)$. Supercuspidal representations always have conductor exponent $c(\pi_v)$ at least 2. The twist $\pi_v \otimes \omega'_v$ of $\pi_v$ by a character $\omega'_v$ of $\mathbb{Q}_v^\times$ is also a supercuspidal representation.

### A.2. Proof of Theorem A.1

A representation $\pi$ of $\mathrm{GL}_2(F)$ (where $F$ is a nonarchimedean local field) of conductor 1, or equivalently of conductor exponent $c(\pi) = 0$, must be a spherical principal series representation $\pi = \omega_1 \boxplus \omega_2$ with both characters $\omega_1, \omega_2$ of $F^\times$ unramified, or equivalently of conductor exponent $c(\omega_1) = c(\omega_2) = 0$. The central character is $\omega_\pi = \omega_1 \omega_2$, which has conductor exponent $c(\omega_1 \omega_2) = 0$; if this central character is trivial, then $\omega_2 = \omega_1^{-1}$. It holds that $\omega_1(x) = |x|_v^{it_v}$ and $\omega_2(x) = |x|_v^{-it_v}$, where $t_v$ is such that the Hecke eigenvalue is $\lambda_\pi(p) = p^{it_v} + p^{-it_v}$. We must think of this as being the local component of $f$ at a prime not dividing the level.

For conductor $p$, or equivalently conductor exponent $c(\pi) = 1$, there are two possibilities. The first is that $\pi$ is a special representation $\omega \mathrm{St}$, where $\omega$ is an unramified character (so that $c(\omega) = 0$). The central character is $\omega_\pi = \omega^2$ that has conductor exponent $c(\omega^2) = 0$. We must think of this as the local component of $f$ at a prime that exactly divides the level and such that the prime does not divide the conductor of the central character of $f$. The other possibility is that $\pi = \omega_1 \boxplus \omega_2$ is a nonspherical principal series representation, where either $c(\omega_1) = 1$ and $c(\omega_2) = 0$ or vice versa. The central character is $\omega_\pi = \omega_1 \omega_2$ that has conductor exponent $c(\omega_1 \omega_2) = 1$. We must think of this as the local component of $f$ at a prime that exactly divides the level and such that the prime also divides the conductor of the central character.

For conductor $p^r$ with $r \geqslant 2$, or equivalently conductor exponent $c(\pi) = r$, there are three possibilities. The first possibility is that $\pi = \omega_1 \boxplus \omega_2$ is a nonspherical principal series representation for which $r = c(\pi) = c(\omega_1) + c(\omega_2) = r$. The central character is $\omega_\pi = \omega_1 \omega_2$ that has conductor exponent at most $r$. The second possibility is that $\pi = \omega \mathrm{St}$ is a special representation, where $\omega$ is a ramified character, so that $c(\omega) \geqslant 1$ and $r = c(\pi) = 2c(\omega)$ (so that $r$ is necessarily even). The central character is $\omega_\pi = \omega^2$ that has conductor exponent at most $r$. The third possibility is that $\pi$ is supercuspidal.

The twist by $\chi$ of a principal series representation $\pi = \omega_1 \boxplus \omega_2$ with central character $\omega_\pi = \omega_1 \omega_2$ is $\pi \otimes \chi = \omega_1 \chi \boxplus \omega_2 \chi$ with central character $\omega_{\pi \otimes \chi} = \omega_1 \omega_2 \chi^2$. The conductor exponent of $\pi \otimes \chi$ becomes $c(\pi \otimes \chi) = c(\omega_1 \chi) + c(\omega_2 \chi)$, and the conductor exponent of $\omega_{\pi \otimes \chi}$ becomes $c(\omega_1 \omega_2 \chi^2)$. The twist of a special representation $\pi = \omega \mathrm{St}$ with central character $\omega_\pi = \omega^2$ is $\pi \otimes \chi = \omega \chi \mathrm{St}$ with central character $\omega_{\pi \otimes \chi} = \omega^2 \chi^2$. The conductor exponent of $\pi \otimes \chi$ becomes $c(\pi \otimes \chi) = \max\{1, 2c(\omega \chi)\}$, and the conductor exponent of $\omega_{\pi \otimes \chi}$ becomes $c(\omega^2 \chi^2)$.

*Proof of Theorem A.1.* The proof uses representation theory, and it suffices to prove Theorem A.1 when $q = p$ is prime. Given a newform $f$ of level $dq$ and central character $\psi^2$, there is a cuspidal automorphic representation $\pi = \pi_\infty \otimes \bigotimes_v \pi_v$ of $\mathrm{GL}_2(\mathbb{A}_\mathbb{Q})$, unique up to isomorphism, with conductor exponent $c(\pi_v)$ at each prime $p$ satisfying $p^{c(\pi_v)} \parallel dq$, whose central character $\omega_\pi$ is the idèlic lift of the primitive character $\chi$ that induces $\psi^2$ and whose global newvector $\xi^\circ$ is the adèlic lift of $f$. If $p \nmid dq$, then $\pi_v$ is spherical and the local conductor exponents ought to satisfy $c(\pi_v) = 0$. All the action happens when $p \mid dq$ for which $\pi_v$ is ramified. The central character $\omega_{\pi_v}$ is the local component $\omega_v^2$ of the idèlic lift of $\chi$, and the local conductor exponent $c(\omega_v)$ is such that $p^{c(\omega_v)} \parallel q$. If $\omega_v$ is nonquadratic, then $c(\omega_{\pi_v}) = c(\omega_v^2) = c(\omega_v)$; however, if $\omega_v$ is quadratic, then $c(\omega_{\pi_v}) = 0$, where $c(\omega_v) = 1$. The accurate quantity of the conductor exponent $c(\pi_v \otimes \chi_v)$ is calculable for given $\pi_v$, where $\chi_v$ is the character corresponding to $\overline{\psi}$. It follows that $\mathrm{cond}(f \otimes \overline{\psi}) \mid q^2$ with equality if $f$ is twist-minimal or $d = 1$ via Lemma 1.4 of [16]. Since $\pi_v$ is a generic irreducible admissible representation of $\mathrm{GL}_2(\mathbb{Q}_v)$, one must go through a case-by-case analysis. Using the above classification, the local conductor exponents can be calculated, although one cannot deal with all the cases in one fell swoop.

1. The case where $\psi$ is quadratic:
   (a) If $d = 1$, then $\pi_v \cong \omega_v \mathrm{St}_v$ for some unramified character $\omega_v$ of $\mathbb{Q}_v^{\times}$ by the above classification. The local component of the twist $f \otimes \overline{\psi}$ is $\pi_v \otimes \chi_v \cong \omega_v \chi_v \mathrm{St}_v$. The conductor exponent of $\pi_v \otimes \chi_v$ thus equals 2, which means that $f \otimes \overline{\psi}$ has conductor $p^2$.
   (b) If $d = q$ is a prime $p$, then $\pi_v$ has conductor exponent $c(\pi_v) = 2$. Since $\omega_{\pi_v}$ is trivial, $\pi_v$ can either be principal series, special or supercuspidal. Humphries [32] discusses this problem in detail and classifies everything in a comprehensive manner. The crucial point is that if $\pi_v = \omega_{p,\mathrm{quad}} \boxplus \omega_{p,\mathrm{quad}}$, where $\omega_{p,\mathrm{quad}}$ denotes the quadratic character of conductor exponent 1, then $\pi_v \otimes \chi_v$ is a spherical principal series representation of conductor exponent 0, so that $f \otimes \overline{\psi}$ has conductor 1. If $\pi_v = \omega_v \mathrm{St}_v$ is special, its twist is a special representation of conductor exponent 1, so that $f \otimes \overline{\psi}$ has conductor $p$. In all other cases, $\pi_v$ is twist-minimal and $\pi_v \otimes \chi_v$ has conductor exponent 2, so that $f \otimes \overline{\psi}$ has conductor $p^2$.
2. The case where $\psi$ is nonquadratic:
   (a) We first handle the case $d = 1$. Since $\psi^2$ is primitive modulo $p$, the corresponding character $\omega_{\pi_v}$ has conductor exponent $c(\omega_{\pi_v}) = 1$. By the above classification, $\pi_v = \omega_{p,1} \boxplus \omega_{p,2}$, where $c(\omega_{p,1}) = 1$, $c(\omega_{p,2}) = 0$ and $\omega_{\pi_v} = \omega_{p,1}\omega_{p,2}$. The twist is $\pi_v \otimes \chi_v = \omega_{p,1}\chi_v \boxplus \omega_{p,2}\chi_v$. Since $\omega_{p,1}\omega_{p,2} = \chi_v^{-2}$ (recall $f$ has central character $\psi^2$), we must have $\omega_{p,1}\chi_v = \omega_{p,2}^{-1}\chi_v^{-1}$, which has conductor exponent 1, as does $\omega_{p,2}\chi_v$. Hence the conductor exponent of $\pi_v \otimes \chi_v$ is $c(\omega_{p,1}\chi_v) + c(\omega_{p,2}\chi_v) = 2$, which means that $f \otimes \overline{\psi}$ has conductor $p^2$.
   (b) If $d = p$, then $\pi_v$ has conductor exponent $c(\pi_v) = 2$. All three possibilities for $\pi_v$ are open, and one can check in each case that twisting leaves its conductor exponent unchanged.

This concludes the proof of Theorem A.1.                                                    □

**Acknowledgments.** The author expresses his gratitude to Eren Mehmet Kıral for supervising this work. He also thanks Valentin Blomer, Jack Buttcane, Paul Nelson and Peter Sarnak for having an interest in earlier versions of this article. Appendix A is attributed to Peter Humphries, and the author acknowledges his diversified comments that improved the readability of this article. Special thanks are owed to Alexandre Perozim de Faveri, Yoichi Motohashi, Maksym Radziwiłł and Will Sawin for their assistance.

**Conflicts of Interest.** None.

**Financial Support.** The author acknowledges the support of the Masason Foundation and the Spirit of Ramanujan STEM Talent Initiative.